\documentclass[11pt,reqno]{amsart}
\usepackage[authoryear,sort]{natbib}
\usepackage{amsmath,amssymb}
\usepackage{dsfont}
\usepackage{xcolor}
\usepackage{graphicx}
\usepackage{bm}
\usepackage[margin=1in]{geometry}
\usepackage[colorlinks=true,citecolor=magenta,linkcolor=blue]{hyperref}
\newcommand{\be}{\begin{equation}}
\newcommand{\ee}{\end{equation}}
\newenvironment{equations}{\equation\aligned}{\endaligned\endequation}
\newcommand{\x}{\mathbf{x}}
\newcommand{\eps}{\varepsilon}
\newcommand{\eq}{\mathrm{eq}}
\newcommand{\init}{\mathrm{in}}

\newcommand{\R}{\mathbb{R}}

\newcommand{\dd}{\mathrm{d}}

\newtheorem{remark}{Remark}

\makeatletter
\let\myorg@bibitem\bibitem
\def\bibitem#1#2\par{%
  \@ifundefined{bibitem@#1}{%
    \myorg@bibitem{#1}#2\par
  }{%
    \begingroup
      \color{\csname bibitem@#1\endcsname}%
      \myorg@bibitem{#1}#2\par
    \endgroup
  }%
}

\makeatother

\makeatletter
\def\blfootnote{\xdef\@thefnmark{}\@footnotetext}
\makeatother

\author[A. Bondesan]{Andrea Bondesan}
\address{Andrea Bondesan \hfill\break
	Department of Mathematical, Physical and Computer Sciences \hfill\break 
    	University of Parma \hfill\break
	Parco Area delle Scienze 53/A, 43124 Parma, Italy \vspace*{2mm}}
 
\email{andrea.bondesan@unipr.it \vspace*{5mm}}

\author[M. Menale]{Marco Menale}
\address{Marco Menale \hfill\break
	Department of Mathematics and Applications ``R. Caccioppoli'' \hfill\break
	University of Naples Federico II \hfill\break
	Via Cintia, Monte S. Angelo, 80126 Napoli, Italy \vspace*{2mm}}
\email{marco.menale@unina.it \vspace*{5mm}}

\author[G. Toscani]{Giuseppe Toscani}
\address{Giuseppe Toscani \hfill\break
	Department of Mathematics ``F. Casorati'' \hfill\break
	University of Pavia \hfill\break
	Via Ferrata 5, 27100 Pavia, Italy \vspace*{2mm}}
\email{giuseppe.toscani@unipv.it \vspace*{5mm}}

\author[M. Zanella]{Mattia Zanella}
\address{Mattia Zanella \hfill\break
	Department of Mathematics ``F. Casorati'' \hfill\break
	University of Pavia \hfill\break
	Via Ferrata 5, 27100 Pavia, Italy \vspace*{2mm}}
\email{mattia.zanella@unipv.it \vspace*{5mm}}

\title[Lotka--Volterra-type kinetic equations]{Lotka--Volterra-type kinetic equations\\for interacting species}

\begin{document}

\vspace*{-0.8cm}
\begin{abstract}
In this work, we examine a kinetic framework for modeling the time evolution of size distribution densities of two populations governed by predator--prey interactions. The model builds upon the classical Boltzmann-type equations, where the dynamics arise from elementary binary interactions between the populations. The model uniquely incorporates a linear redistribution operator to quantify the birth rates in both populations, inspired by wealth redistribution operators. We prove that, under a suitable scaling regime, the Boltzmann formulation transitions to a system of coupled Fokker--Planck-type equations. These equations describe the evolution of the distribution densities and link the macroscopic dynamics of their mean values to a Lotka--Volterra system of ordinary differential equations, with parameters explicitly derived from the microscopic interaction rules. We then determine the local equilibria of the Fokker--Planck system, which are Gamma-type densities, and investigate the problem of relaxation of its solutions toward these kinetic equilibria, in terms of their moments' dynamics. The results establish a bridge between kinetic modeling and classical population dynamics, offering a multiscale perspective on predator--prey systems.
\end{abstract}

\maketitle

\vspace*{0.5cm}
\noindent{\bf Keywords:} Kinetic theory; Lotka--Volterra dynamics; Fokker--Planck equations; Multiscale modeling; Predator--prey interactions.

\vspace*{0.5cm}
\noindent{\bf Mathematics Subject Classification:} 35Q20; 35Q84; 92D25.

%\tableofcontents

%%%%%%%%%%%%%%%%%%%%%%%%%%%%%%%%%%%%%%%%%%%%%%%%%%%%%%%%%%%%%%%%%
%%%%%%%%%%%%%%%%%%%%%%%%%%%%%%%%%%%%%%%%%%%%%%%%%%%%%%%%%%%%%%%%%
%%%%%%%%%%%%%%%%    SECTION 1: INTRODUCTION    %%%%%%%%%%%%%%%%%%
%%%%%%%%%%%%%%%%%%%%%%%%%%%%%%%%%%%%%%%%%%%%%%%%%%%%%%%%%%%%%%%%%
%%%%%%%%%%%%%%%%%%%%%%%%%%%%%%%%%%%%%%%%%%%%%%%%%%%%%%%%%%%%%%%%%

\section{Introduction}
Lotka--Volterra equations, often referred to as predator--prey equations, are first-order nonlinear differential equations that describe population dynamics within systems where multiple species with distinct characteristics interact. In 1925, in the United States, Alfred Lotka proposed a model to describe a chemical reaction with oscillating concentrations. A year later, in Italy, mathematician Vito Volterra independently arrived at the same set of equations while attempting to explain the observed increase in predatory fish, and the corresponding decrease in prey fish, in the Adriatic Sea during World War I. The model proposed by Lotka and Volterra is considered to be the simplest one describing predator--prey interactions for competing species \cite{Hirsch_1988}. A step beyond simple growth dynamics, which are typically influenced by the carrying capacity of the habitat, requires analyzing interacting populations and recognizing their mutual influence on growth. These interactions shape their mutual evolution, making it crucial to predict their dynamics to better understand how communities organize and sustain themselves. 
%In this direction, a rich literature concentrating on species interacting through space-dependent factors have been developed in recent years, see e.g. \cite{BKS,CHS18,DesSor}, in which minimal models for collective behavior may include attraction-repulsion forces between individuals as well as diffusion operators. In fact, cross-interaction can be suitably modeled by nonlinear cross-diffusion factors, which is used to take into account the population pressure \cite{ChoLuiYam_Strong,HHMM,Matthes,Yam}. We also mention selection-type dynamics as generalizations of Lotka--Volterra models \cite{DesJabMisRao,Mirr,Pouchol01012018} which can be seen as a replicator-type dynamics \cite{AFMS} for infinite systems of agents. The study of the evolution of competing species had a tremendous impact in the community studying adaptation dynamics, see e.g. \cite{BERESTYCKI2014264,Berestycki_2009,CalcuaDesRao,Canizo,ChoLuiYam_Weak,DIEKMANN2005257,KISS2008535,Jabin_11}. These approaches have been further generalized in \cite{Lorz17012011,Jabin_2017} in a nonlocal setting. Particle dynamics for cross diffusion systems have been recently introduced in the literature for multicomponent space-dependent interactions \cite{chen21,CDJ,RBM_DAUS}.
In this direction, a rich literature concentrating on species interacting through space-dependent factors have been developed in recent years, see e.g. \cite{BKS,CHS18,DesSor}, in which minimal models for collective behavior may include attraction-repulsion forces between individuals as well as diffusion operators. In fact, cross-interactions can be suitably modeled by nonlinear cross-diffusion factors, which are used to take into account for the population pressure \cite{ChoLuiYam_Strong,HHMM,Matthes,Yam}. We also mention selection-type dynamics as generalizations of Lotka--Volterra models \cite{DesJabMisRao,Mirr,Pouchol01012018} which can be seen as replicator-type dynamics \cite{AFMS} for infinite systems of agents. The study of the evolution of competing species has then tremendously impacted the community investigating adaptation dynamics, see e.g. \cite{BERESTYCKI2014264,Berestycki_2009,CalcuaDesRao,Canizo,ChoLuiYam_Weak,DIEKMANN2005257,KISS2008535,Jabin_11}. These approaches have been further generalized in \cite{Lorz17012011,Jabin_2017} in a nonlocal setting. Moreover, particle dynamics for cross-diffusion systems have been recently introduced in the literature for multicomponent space-dependent interactions \cite{chen21,CDJ,RBM_DAUS}.

Population dynamics have a significant impact in heterogeneous fields, such as ecological systems,  game theory \cite{CRESSMAN2003233} and socio-economic sciences \cite{Palomba-1939}, including for the study of banking systems \cite{SNN-2014}. Additionally, other modeling settings presented in the literature might be considered as mesoscopic approaches inspired by generalized Lotka--Volterra equations \cite{Malcai-2002,Solomon-2002}. In these early approaches, the population size of a system of interacting species is interpreted, within economic systems, with the wealth
of an individual investor or the market value of a traded firm.

Denoting with $X(t) \geq 0$ and $Y(t)\geq 0$ the densities of preys and predators at time $t \geq 0$, the Lotka--Volterra dynamics is classically described by the pair of differential equations
\begin{equations} \label{eq:Lotka-Volterra original}
\frac{\dd X(t)}{\dd t} &= \alpha X(t) - \beta X(t)Y(t), \\[2mm]
\frac{\dd Y(t)}{\dd t} &= -\delta Y(t) + \gamma X(t)Y(t),
\end{equations}
where the evolution of the preys’ population is characterized by a rate of Malthusian growth $\alpha > 0$ and by a predation rate $\beta > 0$, while that of the predators’ population depends on losses due to a death rate $\delta > 0$ and on a growth rate $\gamma > 0$ due to predation.

In this work, we introduce a system of Boltzmann-type kinetic equations to study the time evolution of the size distribution densities of two interacting populations following the rules of a predator–prey system. The dynamics of the introduced model is based on a binary interaction scheme between the introduced populations, the number of births and deaths of individuals in both populations is quantified according to a mass preserving linear operator that has been introduced in a cognate field \cite{BisSpiTos}. In the following, to comply with a statistical description of the dynamics, we will always assume that the number of agents is infinite.  Furthermore, extending the results of \cite{TosZan}, we will show  that, in a quasi-invariant regime of the parameters, corresponding to a high frequency of interactions balanced by a scaling factor, the Boltzmann model reduces to a system of coupled Fokker--Planck-type equations for the distribution densities of the sizes of the populations, where both the diffusion and the drift operators are characterized by time-dependent coefficients. The dynamics of the mean  values of the solutions to this system of Fokker--Planck equations satisfies the Lotka--Volterra model defined in \eqref{eq:Lotka-Volterra original}. Hence, the structure of the latter, along with its relevant parameters, can be justified using a purely microscopic binary interaction scheme, whose aggregate description is formulated assuming the propagation of chaos, through a system of Boltzmann-type equations for the evolution of the particles' densities \cite{CerIllPul}. In particular, this kinetic description allows for the statistical investigation of the underlying multiagent dynamics and the associated Fokker--Planck model provides a simple framework to explicitly determine the relevant kinetic quasi-equilibrium distributions, which take the form of Gamma-type densities, making possible to investigate the convergence of the system's solutions toward equilibrium. In this context, we will study the two significant cases of Malthusian and logistic growths for the preys' population, by deriving the corresponding kinetic models and by analyzing the problem of relaxation to equilibrium in terms of the moments of the solutions. The considered approach is therefore not able to account for the emergence of correlations and possible fluctuations in population sizes near extinction, as shown in \cite{Dobramsyl,Tauber}.

%contains more information on the dynamics of predator--prey interactions than the Lotka--Volterra model \eqref{eq:Lotka-Volterra original}, which appears here as a macroscopic system of equations for the mean values. %Therefore, there is hope that the kinetic systems at hand would furnish important details on how the distributions of predator--prey's populations evolve, including their periodic or limit shapes. 

%Most of these questions were not investigated in \cite{TosZan}, which was mainly devoted to the modeling features behind both the Boltzmann and the Fokker--Planck kinetic systems, in presence of a particular type of elementary binary interactions. More precisely, the choice of these interactions leading to the kinetic system referred only to a situation in which the randomness is proportional to the size of the populations, an assumption which is not always well-justified in the biological setting. For this reason, we shall consider a more general type of elementary interactions  aiming to understand how the choice at the microscopic scale is reflected to the macroscopic behavior of the kinetic systems.

The construction of the general kinetic model of Boltzmann-type will be presented in Section \ref{sec:model}. In Section \ref{sec:FP} we will derive the associated coupled Fokker--Planck system, the analysis of which will form the core of this work. The presence of oscillating time-dependent coefficients in both the diffusion and drift operators, unlike the classical case, prevents the formation of uniform-in-time equilibrium solutions. For this reason, we will concentrate our attention on the role of the so-called quasi-equilibrium distributions, namely the time-dependent distributions that are flux-vanishing for the Fokker--Planck system written in divergence form. These distributions turn out to be Gamma-type densities, whose parameters are functions of the mean values of the solutions to the Fokker--Planck system itself, which in turn solve the Lotka--Volterra model \eqref{eq:Lotka-Volterra original}. By studying the main features of these quasi-equilibrium distributions, we show how they relate to the true solutions of the kinetic system, furnishing a number of interesting information about their evolution, at least for sufficiently large times. Finally, in Section \ref{sec:logistic} we deal with a correction to our kinetic model, switching from a Malthusian-type growth for the preys to a Verhulst-type one. This choice allows to obtain a Fokker--Planck system having still time-dependent drift--diffusion coefficients, but displaying a unique global equilibrium state that represents a large-time attractor for the solutions. 

%\blfootnote{\textit{Date:} \today.}

%%%%%%%%%%%%%%%%%%%%%%%%%%%%%%%%%%%%%%%%%%%%%%%%%%%%%%%%%%%%%%%%%%%
%%%%%%%%%%%%%%%%%%%%%%%%%%%%%%%%%%%%%%%%%%%%%%%%%%%%%%%%%%%%%%%%%%%
%%%%%%%%%%%%%%%    SECTION 2: THE KINETIC MODEL    %%%%%%%%%%%%%%%%
%%%%%%%%%%%%%%%%%%%%%%%%%%%%%%%%%%%%%%%%%%%%%%%%%%%%%%%%%%%%%%%%%%%
%%%%%%%%%%%%%%%%%%%%%%%%%%%%%%%%%%%%%%%%%%%%%%%%%%%%%%%%%%%%%%%%%%%
	
\section{A kinetic model for Lotka--Volterra systems with periodic orbits}\label{sec:model}

We consider an interacting system composed of two populations usually termed as \emph{preys} and \emph{predators}. According to classical kinetic theory, and following the approach developed in \cite{TosZan}, we assume that the relevant microscopic trait characterizing both populations is the size of the interacting groups and that the elementary interactions between groups of preys and predators, which influence the evolution of the overall system, are of binary-type. We suppose in particular that the group sizes can assume all nonnegative real values and we will denote with $x \in \R_+$ and $y \in \R_+$ the number of individuals in each interacting group of preys and predators respectively. For future convenience, we also introduce the vector notation $\x = (x,y)$. Note that, within this assumption, also very small sizes of the populations $x,y \ll 1$ are allowed, a condition which is not fully justified and that has to be suitably avoided. The evolution over time $t \geq 0$ of the two populations of preys and predators can be described by two distribution functions $f_1(x,t)$ and $f_2(y,t)$, depending on their respective sizes $x, y \in \R_+$. More precisely, $f_1(x,t)\dd x$ indicates the fraction of preys belonging to a group of size in $[x, x + \dd x]$ at time $t \geq 0$, while $f_2(y, t)\dd y$ corresponds to the fraction of predators belonging to a group of size in $[y,y+\dd y]$ at time $t \geq 0$. In particular, without loss of generality we can fix the initial distributions $f_1(x,0)$ and $f_2(y,0)$ to be probability density functions, i.e.
\begin{equation*}
\int_{\R_+} f_1(x,0) \dd x = \int_{\R_+} f_2(y,0) \dd y = 1.
\end{equation*}
%
%Then, the related distribution functions are introduced:
%
%\begin{itemize}
    %\item $f_1(t,x)$ is the distribution function, at time $t\geq 0$, of the prey with amount $x$; then, %$f_1(t,x)\, \dd x$ represents the fraction of prey with size in $[x,\, x+\dd x]$, at time $t \geq 0$.
    %\item $f_2(t,y)$ is the distribution function, at time $t\geq 0$, of the predators with amount $y$; %then, $f_2(t,y)\, \dd y$ represents the fraction of predators with size in $[y,\, y+\dd y]$, at time $t \geq %0$.
%\end{itemize}
%Moreover, we assume that the above distribution functions are probabilities at the beginning, i.e. their total mass is equal to $1$. 
The main goal is to characterize the evolution over time of the densities $f_1(x,t)$ and $f_2(y,t)$, from which one can determine the dynamics of their main relevant observable quantities, namely the means and the variances of the sizes of the two species, defined by
\begin{align}
    & m_1(t) = \int_{\R_+} x f_1(x,t) \dd x,
    & m_2(t) = \int_{\R_+} y f_2(y,t) \dd y, \qquad \nonumber \\
    \label{eq:means and energies} \\
    v_1(t) &= \int_{\R_+} (x - m_1(t))^2 f_1(x,t) \dd x,
    & v_2(t) = \int_{\R_+} (y - m_2(t))^2 f_2(y,t) \dd y. \nonumber
\end{align}

In the kinetic approach, the core of the methodology consists in introducing the elementary (binary) interaction rules which are mainly responsible for variations in populations' sizes, in order to write, like in the classical case of gas-dynamics \cite{CerIllPul}, a kinetic system of Boltzmann-type equations modeling the temporal evolution of the distributions of preys and predators, as well as of the associated macroscopic observables defined in \eqref{eq:means and energies}. 

\subsection{Microscopic dynamics}
We start by considering only cross interactions.
These interactions are designed in agreement with the expected behavior of the biological system under consideration: prey--predator interactions should lead to a decrease in the number of preys and, conversely, to an increase in the number of predators, in the presence of a sufficiently large number of preys. Moreover, to take into account possible unpredictable variations, these processes should also be subject to random effects. Following the assumptions in \cite{TosZan}, we define the instantaneous microscopic transitions $x \mapsto x'$ and $y \mapsto y'$ of the form
\begin{equation} \label{eq:microscopic Boltzmann}
    \begin{split}
        x' &= x - \Phi(y)x + \left(\frac{x}{\bar{x}}\right)^p \mathds{1}(x \geq (1-p)s_0) \eta_1(y),
        \\[2mm]
        y' &= y + \Psi(x)y + \left(\frac{y}{\bar{y}}\right)^p \mathds{1}(y \geq (1-p)s_0) \eta_2(x),
    \end{split}
\end{equation}
where $\mathds{1}(A)$ denotes the characteristic function of the set $A \subseteq \R_+$, which is needed to guarantee the positivity of the post-interactions sizes $x'$ and $y'$ (see Remark \ref{remark:positivity of interactions}). Here, the functions $\Phi(y)$ and $\Psi(x)$, which translate deterministic variations within the predator--prey interactions, are assumed as Holling-type II functional responses
\begin{equations} \label{eq:deterministic effects}
        \Phi(y) &= \beta \frac{y}{1+y}, \\[2mm] 
        \Psi(x) &= \gamma \frac{x-\mu}{1+x},
\end{equations}
These functions model respectively the predation rate, i.e. the pressure of predators on the population of preys, and the predators' growth rate due to predation. Here, $\beta, \gamma \in (0,1)$,  and $\mu \geq 1$ are positive constants measuring the maximal predation rate, the maximal growth rate of predators, and the minimal amount of preys that is needed to induce a growth in the predators' population. Note in particular that  the predation and the growth rates remain smaller than one. Moreover, we require that $\gamma$ and $\mu$ satisfy the upper bound $\gamma \mu < 1$, ensuring that even in the absence of preys the number of predators does not collapse to zero in a single interaction. Such type  of predator--prey interactions are in agreement with the pertinent literature on the subject (cf. \cite{holling1959components, holling1966functional, beddington1975mutual, kuno1987principles, abrams2000nature}). Furthermore, relations \eqref{eq:microscopic Boltzmann} also involve stochastic variations through the presence of two independent random variables $\eta_1$ and $\eta_2$, having zero mean and variable bounded variances of the form
\begin{equation} \label{eq:random effects}
    \begin{split}
        \langle\eta_1^2(y)\rangle &= \sigma_1\frac{y}{1+y}, \\[2mm]
        \langle\eta_2^2(y)\rangle &= \sigma_2\frac{x}{1+x},
    \end{split}
\end{equation}
where we have denoted with $\langle \cdot \rangle$ the integration with respect to the probability densities of $\eta_1$ and $\eta_2$, and with $\sigma_1, \sigma_2$ some constant positive parameters. These random effects are activated only when $x$ and $y$ are greater than $(1-p) s_0$, with $0 < s_0 < 1$, and encode the  diffusion of preys and predators, which is assumed to depend in a nonlinear way from the presence of the other species. One can think for example of preys diffusing as a way to escape from or avoid predators and, conversely, predators diffusing to reach more preys. The choice once again of a Holling-type II form for the variances of $\eta_1$ and $\eta_2$ translates the fact that the ability of one species to diffuse decreases with the increase of the other species' population, up to a maximum value prescribed by the constants $\sigma_1$ and $\sigma_2$. The nonlinearity of such random effects is further modulated by the parameter $p$, satisfying $0 < p \leq 1$ (note that in the work \cite{TosZan}, only the case $p=1$ has been considered), which tells how strong the species' diffusion is, based on the reference values $\bar{x}$ for the preys and $\bar{y}$ for the predators. In particular, for $p < 1$, these thresholds define positive bounds on the species' sizes below which the random effects are more evident, while they are mitigated when $x > \bar{x}$ and $y > \bar{y}$. Additionally, as already done for existing models with similar formulations \cite{ParTos-2013}, one can verify that the choice $p < 1$ leads to formation of long-time equilibrium distributions with thin tails, while the choice $p = 1$ produces fat-tails equilibrium distributions. For this reason, and without loss of essential generality, in what follows we will consider either $p = 1/2$ or $p = 1$ to quantify the prototypical behavior of the tails that characterize the local equilibria of the kinetic system. Moreover, in order to ease the notations in the remaining part, we shall assume from now on that $\bar{x} = \bar{y} = 1$. 

\begin{remark} \label{remark:positivity of interactions}
We highlight that the microscopic transitions \eqref{eq:microscopic Boltzmann} must map positive values into positive values, while the random variables $\eta_1(y)$ and $\eta_2(x)$ may have both positive and negative values, since they are chosen to have zero mean. Clearly, the positivity of the updates $x'$ and $y'$ is always guaranteed in the case $p = 1$ considered in \cite{TosZan}, where the indicator functions are not playing a role. Indeed, if $\eta_1(y) $ satisfies the negative lower bound $\eta_1(y) \geq - (1-\beta) \bar{x}$, then
\begin{equation*}
    x' \geq x \left(1 - \beta + \frac{1}{\bar{x}} \eta_1(y)\right) \geq 0.
\end{equation*}
Similarly, the positivity of $y'$ would follow by assuming that the random variable $\eta_2(x)$ satisfies the lower bound $\eta_2(x) \geq - (1 - \gamma\mu) \bar{y}$. On the contrary, the use of the indicator functions inside the interactions \eqref{eq:microscopic Boltzmann} cannot be avoided when $p < 1$, since otherwise the possibility for $\eta_1(y)$ and $\eta_2(x)$ to be negative in an interval of strictly positive measure would not allowed. Indeed, in absence of the indicator function, for a given $p < 1$, one would compute
\begin{equation*}
    x' \geq x \left(1 - \beta + \frac{1}{x^{1-p} {\bar{x}}^p} \eta_1(y)\right),
\end{equation*}
and infer the $x$-dependent bound
\begin{equation*}
    \eta_1(y) \geq - (1-\beta) x^{1-p} {\bar{x}}^p.
\end{equation*}
This unpleasant fact can be circumvented as soon as $x \geq (1-p)s_0$, since in this case $\eta_1(y)$ would satisfy the uniform lower bound
\begin{equation*}
    \eta_1(y) \geq - (1-\beta)[(1-p)s_0]^{1-p} {\bar{x}}^p,
\end{equation*}
whence the introduction of the indicator functions in \eqref{eq:microscopic Boltzmann} to cut off all random effects below the positive threshold $(1-p)s_0$. Notice that the latter has been chosen of this particular form simply to highlight that the restriction is no longer present when $p = 1$, while the coefficient $s_0 < 1$ is there to ensure that we can reduce the threshold at will.
\end{remark}

Next, to take into account expected growth and shrinking factors of the two populations, we proceed as in \cite{BISI_17,BisSpiTos} and introduce a redistribution-type operator encapsulating birth and death processes. The latter processes can be linked to interactions of the two species with an additional surrounding background, that can be thought of as the resources (for example food, shelters, etc.) offered by the environment. If we denote with $z \in \R_+$ the variable representing the size of these resources, and with $g(z,t)$ and $h(z,t)$ their corresponding kinetic densities, the sought redistribution operator acts at the particle level as
\begin{equation} \label{eq:microscopic redistribution}
    \begin{split}
        x'' &= x + \alpha \left(z - \chi x \right), \\[2mm]
        y'' &= y + \nu \left(z - \theta y \right),
    \end{split}
\end{equation}
where $\alpha, \nu > 0$ are constant interaction coefficients, while $\chi, \theta > -1$ define some constant rates of variation for the two populations. In particular, we require that these four parameters satisfy the constraints $\alpha\chi < 1$ and $\nu\theta < 1$, in order to ensure that relations \eqref{eq:microscopic redistribution} preserve the positivity of the group sizes. Furthermore, we assume that the means of the densities $g(z,t)$ and $h(z,t)$ are explicitly given by
\begin{equation}
\label{eq:mean of environment}
\int_{\mathbb R_+} z g(z,t) \dd z = (\chi+1) m_1(t), \qquad \int_{\mathbb R_+} z h(z,t) \dd z = (\theta+1) m_2(t).
\end{equation}
We highlight that since $\alpha$ and $\nu$ are assumed to be positive, the presence of the environment always ensures some growth (even if, possibly, very small) of the populations of preys and predators, alike. As we shall see, this property is essential to recover well-behaved equilibrium states for the kinetic model. From a modeling perspective, we interpret this assumption by saying that preys tend to reproduce by having access to food or by using shelters to avoid the predators, while the latter may reproduce by feeding on other types of preys or (when handling, for example) by escaping from predation themselves, taking advantage of the environment to hide. Additionally, the assumption on the means of $g$ and $h$ convey the idea that the environment has not enough ($-1 < \chi,\theta < 0$), enough ($\chi = \theta = 0$) or more than enough ($\chi,\theta > 0$) resources to cover the needs, of food or space, of each population.

%%%%%%%%%%%%%%%%%%%%%%%%%%%%%%%%%%%%%%%%%%%%%%%%%%%%%%%%%%%%%%%%%%%
\subsection{The kinetic model}
Being interested in a statistical description of the dynamics, in the following we will always assume that the number of agents/particles is infinite. The analysis that follows is conducted under the ansatz of the propagation of chaos  \cite{CerIllPul}.
Starting from the interaction rules \eqref{eq:microscopic Boltzmann}, we define the system of kinetic equations modeling the evolution of the two distribution functions, $f_1(x,t)$ and $f_2(y,t)$. This evolution is provided by the following Boltzmann-type equations \cite{ParTos-2013, TosZan}
\begin{equation} \label{eq:Boltzmann}
    \begin{split}
        \displaystyle \frac{\partial f_1(x,t)}{\partial t} = R_{\chi}^{\alpha}(f_1)(x,t)+Q_{12}(f_1,f_2)(x,t), \\[2mm] 
        \displaystyle \frac{\partial f_2(y,t)}{\partial t} = R_{\theta}^{\nu}(f_2)(y,t)+Q_{21}(f_2,f_1)(y,t),
    \end{split}
\end{equation}
where $R_\chi^\alpha(f_1)$ and $R_\theta^\nu(f_2)$ are redistribution-type linear operators describing the balance of the microscopic interactions \eqref{eq:microscopic redistribution}, whereas $Q_{12}(f_1,f_2)$ and $Q_{21}(f_2,f_1)$ are Boltzmann-type bilinear operators accounting for the microscopic interactions defined by \eqref{eq:microscopic Boltzmann}. The Boltzmann system \eqref{eq:Boltzmann} can be conveniently recast in weak form, which corresponds to evaluate the action that the solutions to system \eqref{eq:Boltzmann}  exercises on the observable quantities, which are in general characterized by two  smooth test functions $\varphi(x)$ and $\psi(y)$. This leads to the system
\begin{equation} \label{eq:Boltzmann weak}
    \begin{split}
        \frac{\dd}{\dd t}\int_{\R_+ }\varphi(x)f_1(x,t) \dd x &= \int_{\R_+ \times \R_+} \big(\varphi(x'')-\varphi(x)\big) f_1(x,t) g(z,t) \dd z \dd x \\[2mm]
        & \hspace*{1cm} +\int_{\R_+\times \R_+} \kappa_1(x,y) \big\langle\varphi(x')-\varphi(x)\big\rangle f_1(x,t)f_2(y,t) \dd x \dd y, \\[4mm]
        \frac{\dd}{\dd t}\int_{\R_+}\psi(y)f_2(y,t) \dd y 
        &= \int_{\R_+\times \R_+} \big(\psi(y'')-\psi(y)\big) f_2(y,t) h(z,t) \dd z \dd y \\[2mm]
        & \hspace*{1cm} +\int_{\R_+\times\R_+} \kappa_2(x,y) \big\langle\psi(y')-\psi(y)\big\rangle f_1(x,t)f_2(y,t) \dd x \dd y.
    \end{split}
\end{equation}
In system \eqref{eq:Boltzmann weak}, $\kappa_i(x,y) > 0$, $i = 1,2$ denote the interaction frequencies that prescribe the rate at which the two species encounter each other. They are assumed of the following form
\begin{equation} \label{eq:interaction frequencies}
    \begin{split}
        \kappa_1(x,y) &= \kappa(y) = 1+y, \\[2mm]
        \kappa_2(x,y) &= \kappa(x) = 1+x.
    \end{split}
\end{equation}
By this choice, in the spirit of the Lotka--Volterra model, the frequency of interactions in the first equation grows with the number of predators, while in the second equation it grows with the number of preys.

%%%%%%%%%%%%%%%%%%%%%%%%%%%%%%%%%%%%%%%%%%%%%%%%%%%%%%%%%%%%%%%%%%%
\subsection{Evolution of the observable quantities}

In particular, if one assumes $\varphi(x) = x^r$ and $\psi(y) = y^r$, for $r \in \mathbb{N}$, the kinetic equations \eqref{eq:Boltzmann weak} provide the time evolution of the $r$th-order moment of the system. For the aim of this paper, we will focus on the $1$st-order moment and the $2$nd-order moments, defined for each of the two species by \eqref{eq:means and energies}. The former represents the mean, while the latter provides the energy (and thus also the variance) of the system.

At first, we observe that taking $\varphi(x) = 1$ and $\psi(y) = 1$, one deduces from equations \eqref{eq:Boltzmann weak} the evolution of the densities of the two populations, i.e.
\begin{equation*}
    \frac{\dd}{\dd t}\int_{\R_+} f_1(x,t) \dd x = 
    \frac{\dd}{\dd t}\int_{\R_+} f_2(y,t) \dd y = 0,
\end{equation*}
ensuring the conservation of the total mass of the system. Since $f_1(x,t)$ and $f_2(y,t)$ can be chosen to be probability densities at initial time, this conservation property then ensures that, for all $t \geq 0$,
\begin{equation} \label{eq:conservation of mass}
    \int_{\R_+} f_1(x,t) \dd x = \int_{\R_+} f_2(y,t) \dd y = 1,
\end{equation}
meaning that the two distributions remain probability densities during the overall evolution of the system. Indeed, notice that the positivity of solutions to system \eqref{eq:Boltzmann} can be readily proved by resorting to standard arguments developed for the classical Boltzmann equation \cite{CerIllPul}, owing to the fact that the redistribution operators $R_\chi^\alpha$ and $R_\theta^\nu$ are positivity-preserving transport-like operators and to the specific form of the interaction frequencies \eqref{eq:interaction frequencies}. This choice allows to split the Boltzmann-like operators $Q_{12}$ and $Q_{21}$ into a gain and a loss part, where the former is positivity-preserving while the latter is a multiplicative-like operator acting linearly on $f_1$ and $f_2$, and depends on the other species' distribution only through its zeroth- and first-order moments, via integration in $y$ and in $x$, respectively, of the quantities $\kappa(y) f_2(y,t)$ and $\kappa(x) f_1(x,t)$. In particular, these moments are positive since from \eqref{eq:conservation of mass} the masses of the populations are conserved and their means solve the Lotka--Volterra system \eqref{eq:Lotka-Volterra}, as we shall observe in what follows.

If now $\varphi(x) = x$ and $\psi(y) = y$, from \eqref{eq:Boltzmann weak} we obtain the coupled evolution of the mean values $m_1(t)$ and $m_2(t)$ defined in \eqref{eq:means and energies}. In particular, we get
\begin{equation}\label{eq:Lotka-Volterra}
\begin{split}
\frac{\dd m_1(t)}{\dd t} &= \alpha m_1(t) - \beta m_1(t) m_2(t), \\[2mm]
\frac{\dd m_2(t)}{\dd t} &= - (\gamma\mu - \nu) m_2(t) + \gamma m_1(t) m_2(t), 
\end{split}
\end{equation}
which, under the identification $\delta = \gamma\mu - \nu$ and assuming that $\nu < \gamma\mu$, corresponds to the classical Lotka--Volterra model \eqref{eq:Lotka-Volterra original}. For any initial data $m_1(0)$, $m_2(0) > 0$, a unique vector solution $\mathbf{m}(t) = (m_1(t), m_2(t))$ exists for all times and remains positive (see for example \cite{IanPug}). The solution $\mathbf{m}(t)$ evolves, as known,  in the phase plane, over closed orbits around the unique non-zero equilibrium point 
\begin{equation} \label{eq:m*}
    \mathbf{m}^* = \left( \frac{\delta}{\gamma}, \frac{\alpha}{\beta} \right).
\end{equation}
%In particular, if we denote with $T > 0$ the period of a solution to \eqref{eq:Lotka-Volterra}, it is possible to show \cite{IanPug} that the averages of $m_1$ and $m_2$ over one period are conserved and coincide with the corresponding values of the equilibrium point $\mathbf{m}^\infty$, namely
%\begin{align*}
%    \overline{m}_1(t) := \int_0^T m_1(t) \dd t = \frac{\delta}{\gamma}, \\[2mm]
%    \overline{m}_2(t) := \int_0^T m_2(t) \dd t = \frac{\alpha}{\beta}.
%\end{align*}
In particular, the function $H(m_1(t),m_2(t)) = \gamma m_1(t) - \delta \log m_1(t) + \beta m_2(t) - \alpha \log m_2(t)$ is a conserved quantity of the system, that characterizes the shape of such orbits spanned by the preys and the predators in the phase plane. The solutions $m_1(t)$ and $m_2(t)$ are thus periodic and bounded from above and below. For future use, we define these bounds via some positive constants $\underline{c}_1$, $\overline{c}_1$ and $\underline{c}_2$, $\overline{c}_2$ such that for any $t \in \R_+$
\begin{equation} \label{eq:bounds on means}
\begin{split}
    & \underline{c}_1 \leq m_1(t) \leq \overline{c}_1, \\[2mm]
    & \underline{c}_2 \leq m_2(t) \leq \overline{c}_2.
\end{split}
\end{equation}

We conclude by studying the evolution of the variances \eqref{eq:means and energies} of $f_1(x,t)$ and $f_2(y,t)$, which are obtained by inserting $\varphi(x) = (x - m_1(t))^2$ and $\psi(y) = (y - m_2(t))^2$ in \eqref{eq:Boltzmann weak}. While the means of the two distributions characterize the orbits of the Lotka--Volterra equations \eqref{eq:Lotka-Volterra}, their variances account for the deviations around these orbits. We distinguish the case $p = 1$, for which we obtain 
%\begin{equation} \label{eq:energy Boltzmann p=1}
%\begin{split}
%\frac{\dd E_1(t)}{\dd t} &= - \big((2\beta - \sigma_1) m_2(t) + 2\alpha\chi \big) E_1(t) + \alpha^2 \chi^2 E_2(t) + 2\alpha (\chi+1) m_1^2(t) \\[2mm]
%& \hspace*{1cm} -2\alpha^2 \chi (\chi+1) m_1^2(t) + \alpha^2 \int_{\R_+} z^2 g(z,t) \dd z + \beta^2 E_1(t) \int_{\R_+} \frac{y^2}{1+y} f_2(y,t) \dd y, \\[4mm]
%\frac{\dd E_2(t)}{\dd t} &= - \big(2\gamma(\mu - m_1(t)) - \sigma_2 m_1(t) + 2\nu \theta \big) E_2(t) + \nu^2 \theta^2 E_2(t) + 2\nu (\theta+1) m_2^2(t) \\[2mm]
%& \hspace*{1cm} -2\nu^2 (\theta+1) m_2^2(t) + \nu^2 \int_{\R_+} z^2 h(z,t) \dd z + \gamma^2 E_2(t) \int_{\R_+} \frac{(x-\mu)^2}{1+x} f_1(x,t) \dd x,
%\end{split}
%\end{equation}
%from the case $p = 1/2$ that gives
%\begin{equation} \label{eq:energy Boltzmann p=1/2}
%\begin{split}
%\frac{\dd E_1(t)}{\dd t} &= -2(\beta m_2(t) + \alpha \chi) E_1(t) + \alpha^2 \chi^2 E_1(t) + 2\alpha (\chi+1) m_1^2(t) - 2\alpha^2 \chi (\chi+1) m_1^2(t) \\[2mm]
%& \hspace*{1cm} + \sigma_1 m_1(t) m_2(t) + \alpha^2 \int_{\R_+} z^2 g(z,t) \dd z + \beta^2 E_1(t) \int_{\R_+} \frac{y^2}{1+y}f_2(y,t) \dd y \\[4mm]
%\frac{\dd E_2(t)}{\dd t} &= -2 \big(\gamma(\mu - m_1(t)) + \nu\theta \big) E_2(t) + \nu^2 \theta^2 E_2(t) + 2\nu (\theta+1) m_2^2(t) - 2\nu^2 (\theta+1) m_2^2(t) \\[2mm]
%& \hspace*{1cm} + \sigma_2 m_1(t) m_2(t) + \nu^2 \int_{\R_+} z^2 h(z) \dd z + \gamma^2 E_2(t) \int_{\R_+} \frac{(x-\mu)^2}{1+x} f_1(x,t) \dd x.
%\end{split}
%\end{equation}
\begin{equation} \label{eq:variance Boltzmann p=1}
    \begin{split}
        \frac{\dd v_1(t)}{\dd t} &= -[(2\beta - \sigma_1) m_2(t) + 2\alpha \chi] v_1(t) + \sigma_1 m_1^2(t) m_2(t) + \alpha^2 \chi^2 (v_1(t) + m_1^2(t)) \\[2mm]
        & \hspace*{1cm} -2 \alpha^2 \chi(\chi+1) m_1^2(t) + \alpha^2 \int_{\R_+} z^2 g(z,t) \dd z + \beta^2 (v_1(t) + m_1^2(t)) \int_{\R_+} \frac{y^2}{1+y} f_2(y,t) \dd y,
        \\[4mm]
        \frac{\dd v_2(t)}{\dd t} &= -[2\gamma(\mu - m_1(t)) - \sigma_2 m_1(t) + 2\nu \theta] v_2(t) + \sigma_2 m_1(t) m_2^2(t) + \nu^2 \theta^2 (v_2(t) + m_2^2(t)) \\[2mm]
        & \hspace*{1cm} -2\nu^2 (\theta+1) m_2^2(t) + \nu^2 \int_{\R_+} z^2 h(z,t) \dd z + \gamma^2 (v_2(t) + m_2^2(t)) \int_{\R_+} \frac{(x - \mu)^2}{1+x} f_1(x,t) \dd x,
    \end{split}    
\end{equation}
from the case $p = 1/2$ that gives
\begin{equation} \label{eq:variance Boltzmann p=1/2}
    \begin{split}
        \frac{\dd v_1(t)}{\dd t} &= -2 (\beta m_2 + \alpha \chi) v_1(t) + \sigma_1 m_1(t) m_2(t) + \alpha^2 \chi^2 (v_1(t) + m_1^2(t)) \\[2mm]
        & \hspace*{1.5cm} -2\alpha^2 \chi (\chi+1) m_1^2(t) + \alpha^2 \int_{\R_+} z^2 g(z,t) \dd z  \\[2mm]
         & \hspace*{2.5cm} +\beta^2 (v_1(t) + m_1^2(t)) \int_{\R_+} \frac{y^2}{1+y} f_2(y,t) \dd y + \int_0^{s_0/2} x f_1(x,t)dx, \\[4mm] 
        \frac{\dd v_2(t)}{\dd t} &= -2 [\gamma(\mu - m_1(t)) + \nu \theta] v_2(t) + \sigma_2 m_1(t) m_2(t) +\nu^2 \theta^2 (v_2(t) + m_2^2(t)) \\[2mm]
        & \hspace*{1.5cm} -2 \nu^2 \theta(\theta+1) m_2^2(t) + \nu^2 \int_{\R_+} z^2 h(z,t) \dd z \\[2mm]
        &\hspace*{2.5cm} + \gamma^2 (v_2(t) + m_2^2(t)) \int_{\R_+} \frac{(x - \mu)^2}{1+x} f_1(x,t) \dd x + \int_0^{s_0/2} y f_2(y,t)dy.
    \end{split}    
\end{equation}
%Therefore, the evolution of each variance, as the distributions $f_1(t,x)$ and $f_2(t,y)$ are probability denties for all $t \geq 0$, is obtained for $\varphi(x)=(x-m_1(t))^2$ and $\varphi(y)=(y-m_2(t))^2$ in the equation \eqref{eqboltz5} and \eqref{eqboltz6}, respectively.
It is worth pointing out that the integrals $\displaystyle \int_{\R_+} z^2 g(z,t) \dd z$ and $\displaystyle \int_{\R_+} z^2 h(z,t) \dd z$ are the $2$nd-order moments of the densities $g(z,t)$ and $h(z,t)$ characterizing the evolution of the background environment, whose means have been assumed to have the explicit form \eqref{eq:mean of environment}. Therefore, these integrals can be computed once also the energies of $g(z,t)$ and $h(z,t)$ are explicitly assigned. For instance, if the densities were taken to be Dirac delta functions, namely a situation where the resources provided by the environment have zero variance and are thus kept fixed, then, according to hypotheses \eqref{eq:mean of environment} one would have
\begin{align*}
    \int_{\R_+} z^2 g(z,t) \dd z & = \left[(\chi+1) m_1(t) \right]^2, \\[2mm]
    \int_{\R_+} z^2 h(z,t) \dd z &= \left[(\theta+1) m_2(t) \right]^2.
\end{align*}
%Bearing the above systems for the energies in mind, the evolution of the related variances is easily recovered from the relations $v_i(t) = E_i(t) - m_i^2(t)$ for $i = 1, 2$. Therefore, for $p=1$ one finds
Therefore, unlike what happens for the evolution of the means $m_1(t)$ and $m_2(t)$, we see from systems \eqref{eq:variance Boltzmann p=1} and \eqref{eq:variance Boltzmann p=1/2} that it is hard to recover a closed evolution for the variances of the species' distributions (nor for their higher order moments). Hence, to simplify the study of the behavior of the solutions to the Boltzmann-type system \eqref{eq:Boltzmann}, in the next section, we propose a reduced complexity approach by resorting to a simplified model of Fokker--Planck-type which allows for an explicit closed evolution of the principal moments, as well as of an explicit evaluation of the local equilibria associated with the kinetic system.

%\begin{align*}
%    &\int_{\R_+}\frac{y^2}{1+y}g(y,t) \dd y \\
%    &\int_{\mathbf{R}^+}\frac{(x-\mu)^2}{1+x}f(x,t) \dd x,
%\end{align*}
%of the equations \eqref{varf}-\eqref{varg}, are not related to any moment of the system.

%%%%%%%%%%%%%%%%%%%%%%%%%%%%%%%%%%%%%%%%%%%%%%%%%%%%%%%%%%%%%%%%%%%
%%%%%%%%%%%%%%%%%%%%%%%%%%%%%%%%%%%%%%%%%%%%%%%%%%%%%%%%%%%%%%%%%%%
%%%%%%%%%%%%%%%%%    SECTION 3: FOKKER-PLANCK    %%%%%%%%%%%%%%%%%%
%%%%%%%%%%%%%%%%%%%%%%%%%%%%%%%%%%%%%%%%%%%%%%%%%%%%%%%%%%%%%%%%%%%
%%%%%%%%%%%%%%%%%%%%%%%%%%%%%%%%%%%%%%%%%%%%%%%%%%%%%%%%%%%%%%%%%%%

\section{Derivation of a coupled Fokker--Planck system}\label{sec:FP}

%%%%%%%%%%%%%%%%%%%%%%%%%%%%%%%%%%%%%%%%%%%%%%%%%%%%%%%%%%%%%%%%%%%
\subsection{Scaling of the redistribution operators}

In this section,  we show that for a certain range of the parameters characterizing the kinetic system \eqref{eq:Boltzmann weak}, usually referred to as the \emph{quasi-invariant collision regime}, it is possible to derive a reduced model of Fokker--Planck type \cite{ParTos-2013}. The Fokker--Planck description in many cases facilitates the analytical study of the shape of the emerging equilibrium distributions. We refer for example to \cite{DPTZ, FPTT19, MR3597010} for further details on this topic.

Let us first concentrate on the redistribution operators $R_\chi^\alpha(f_1)$ and $R_\theta^\nu(f_2)$, resulting from the particle dynamics \eqref{eq:microscopic redistribution}. A particle approach to the redistribution dynamics has been studied in \cite{BISI_17}. By denoting with $\eps \ll 1$ a suitable scaling parameter, we may introduce a new timescale of observation $\tau = \eps t$, together with the rescaled distributions $f_1^\eps(x,\tau) = f_1(x,\tau/\eps)$ and $f_2^\eps(x,\tau/\eps) = f_2(x,\tau/\eps)$ for the two populations, as well as $g^\eps(z,\tau) = g(z,\tau/\eps)$ and $h^\eps(z,\tau) = h(z,\tau/\eps)$ for those characterizing the environment, obtaining the temporal evolution
\begin{equation} \label{eq:rescaled redistribution 1}
    \begin{split}
        \frac{\dd}{\dd \tau} \int_{\R_+} \varphi(x) f^\eps_1(x,\tau) \dd x &= \frac{1}{\eps} \int_{\R_+ \times \R_+} \big( \varphi(x'') - \varphi(x) \big) f_1^\eps(x,\tau) g^\eps(z,\tau) \dd z \dd x, \\[2mm]
        \frac{\dd}{\dd \tau} \int_{\R_+} \psi(y) f^\eps_2(y,\tau) \dd y &= \frac{1}{\eps} \int_{\R_+ \times \R_+} \big( \psi(y'') - \psi(y) \big) f_2^\eps(y,\tau) h^\eps(z,\tau) \dd z \dd y.
    \end{split}
\end{equation}
Hence, considering the following scaling of the parameters
\begin{equation*}
    \alpha \to \eps \alpha, \qquad \nu \to \eps \nu, 
\end{equation*}
the interactions \eqref{eq:microscopic redistribution} are quasi-invariant and we may thus perform the Taylor expansion
\begin{equation*}
    \begin{split}
        \varphi(x'') - \varphi(x) = (x'' - x) \varphi'(x) + \frac{1}{2} (x'' - x)^2 \varphi''(\hat{x}), \\[2mm]
        \psi(y'') - \psi(y) = (y'' - y) \psi'(y) + \frac{1}{2} (y'' - y)^2 \psi''(\hat{y}),
    \end{split}
\end{equation*}
with $\hat{x} \in (\min \{ x, x'' \}, \max \{ x, x'' \})$ and $\hat{y} \in (\min \{ y, y'' \}, \max \{ y, y'' \})$. Plugging the above expansions into \eqref{eq:rescaled redistribution 1}, we then obtain, for any sufficiently smooth test functions $\varphi(x)$ and $\psi(y)$,
\begin{equation*}
    \begin{split}
        \frac{\dd}{\dd \tau} \int_{\R_+} \varphi(x) f_1^\eps(x,\tau) \dd x  &= \int_{\R_+ \times \R_+} \alpha (z - \chi x)\varphi'(x) f_1^\eps(x,\tau) g^\eps(z,\tau) \dd z \dd x + \mathcal{R}_\varphi(f_1^\eps, g^\eps)(x,\tau), \\[2mm]
        \frac{\dd}{\dd \tau} \int_{\R_+} \psi(y) f_2^\eps(y,\tau) \dd y  &= \int_{\R_+ \times \R_+} \nu (z - \theta y) \psi'(y) f^\eps_2(y,\tau) h^\eps(z,\tau) \dd z \dd y + \mathcal{R}_\psi(f_2^\eps, h^\eps)(y,\tau),
    \end{split}
\end{equation*}
with the integral remainders $\mathcal{R}_\varphi$ and $\mathcal{R}_\psi$ given by
\begin{equation*}
    \begin{split}
        \mathcal{R}_\varphi(f_1^\eps, g^\eps)(x,\tau) &= \frac{\eps}{2} \int_{\R_+ \times \R_+} \alpha^2 (z - \chi x)^2 \varphi''(\hat{x}) f_1^\eps(x,\tau) g^\eps(z,\tau) \dd z \dd x, \\[2mm]
        \mathcal{R}_\psi(f_2^\eps, h^\eps)(y,\tau) &= \frac{\eps}{2} \int_{\R_+ \times \R_+} \nu^2 (z - \theta y)^2 \psi''(\hat{y}) f_2^\eps(y,\tau) h^\eps(z,\tau) \dd z \dd y. 
    \end{split}
\end{equation*}
Therefore, observing that
\begin{equation*}
    \left|\mathcal{R}_\varphi(f_1^\eps, g^\eps) \right| \underset{\eps \to 0}{\longrightarrow} 0, \qquad \left|\mathcal{R}_\psi(f_2^\eps, h^\eps) \right| \underset{\eps \to 0}{\longrightarrow} 0,
\end{equation*}
in the limit $\eps \to 0$ the distributions $f^\eps_1(x,\tau)$ and $f^\eps_2(y,\tau)$ formally converge to some $f_1(x,\tau)$ and $f_2(y,\tau)$, which are solutions to the system
\begin{equation*}
    \begin{split}
        \frac{\dd}{\dd \tau}\int_{\R_+} \varphi(x)f_1(x,\tau) \dd x &= \int_{\R_+ \times \R_+} \alpha(z - \chi x) \varphi'(x) f_1(x,\tau) g(z,\tau) \dd z \dd x, \\[2mm]
        \frac{\dd}{\dd \tau} \int_{\R_+} \psi(y) f_2(y,\tau) \dd y  &= \int_{\R_+ \times \R_+} \nu(z - \theta y) \psi'(y) f_2(y,\tau) h(z,\tau) \dd z \dd y.
    \end{split}
\end{equation*}
Integrating back by parts the resulting equations and recalling \eqref{eq:mean of environment}, we recover the following transport equations
\begin{equation} \label{eq:FP redistribution}
    \begin{split}
        \frac{\partial f_1(x,\tau)}{\partial \tau} &= \alpha \frac{\partial}{\partial x} \big((\chi x - (\chi+1) m_1(\tau)) f_1(x,\tau) \big) = \bar{R}_\chi^\alpha(f_1)(x,\tau), \\[2mm]
        \frac{\partial f_2(y,\tau)}{\partial \tau} &= \nu \frac{\partial}{\partial y} \big((\theta y - (\theta+1) m_2(\tau)) f_2(y,\tau) \big) = \bar{R}_\theta^\nu(f_2)(y,\tau), 
    \end{split}
\end{equation}
coupled with the boundary conditions 
\[
\big(\chi x - (\chi+1)m_1(\tau) \big) f_1(x,\tau)\Big|_{x = 0, +\infty} = 0, \quad \big( \theta y - (\theta+1)m_2(\tau) \big) f_2(y,\tau) \Big|_{y = 0, +\infty} = 0.
\]
We highlight that $\bar{R}_\chi^\alpha$ and $\bar{R}_\theta^\nu$ correspond exactly to the redistribution operators introduced in \cite{BisSpiTos}. In particular, in view of the above no-flux boundary conditions, one finds
\begin{align*}
    &\int_{\R_+} \bar{R}_{\chi}^{\alpha}(f_1)(x,\tau) \dd x = 0,
    &&\int_{\R_+} \bar{R}_{\theta}^{\nu}(f_2)(y,\tau) \dd y = 0,
    \\[2mm]
    &\int_{\R_+} x \bar{R}_{\chi}^{\alpha}(f_1)(x,\tau) \dd x = \alpha m_1(\tau),
    &&\int_{\R_+} y \bar{R}_{\theta}^{\nu}(f_2)(y,\tau) \dd y = \nu m_2(\tau), 
\end{align*}
meaning that the redistribution operators preserve the total masses of the two populations, while increasing the mean sizes of the preys by a factor $\alpha > 0$ and of the predators by a factor $\nu > 0$. The new parameters $\chi, \theta > -1$ define the nature of the redistribution. In particular, for positive values of these parameters, the redistribution operators promote the growth of small populations and simultaneously penalize very large populations,
%Conversely, for $\chi, \theta < -1$, the opposite occurs: smaller populations supply additional resources to the larger ones. 
while for $-1 < \chi, \theta < 0$, the redistribution mechanism favors both large and small populations at the expense of medium-sized ones. We can obtain a further insight on this operator by solving the homogeneous equations in \eqref{eq:FP redistribution}. Starting from the initial conditions $f_1(x,0) = f_1^\init(x)$ and $f_2(y,0) = f_2^\init(y)$, with means $m_1(0) > 0$ and $m_2(0) > 0$,
%the redistribution operators naturally increase the mean size of the populations since from \eqref{eq:FP redistribution} we have $\frac{d}{d\tau}{m_1} = \alpha m_1(\tau)$, $\frac{d}{d\tau}{m_2} = \nu m_2(\tau)$  and therefore $m_1(\tau) = m_1(0)e^{\alpha \tau}$, $m_2(\tau) = m_2(0)e^{\nu \tau}$ $\tau\ge0$,  and
the solutions to \eqref{eq:FP redistribution} read
\begin{equation*}
    \begin{split}
        f_1(x,\tau) &= e^{\alpha\chi \tau} f_1^\init \left(e^{\alpha\chi \tau} x - m_1(0) (e^{\alpha(\chi+1) \tau} - 1) \right), \\[2mm]
        f_2(y,\tau) &= e^{\nu\theta \tau} f_2^\init \left(e^{\nu\theta \tau} y - m_2(0) (e^{\nu(\theta+1) \tau} - 1) \right). 
    \end{split}
\end{equation*}
Hence, if $f_1^\init(x) \equiv 0$, for any $x < 0$ the corresponding solution $f_1(x,\tau)$ satisfies 
\begin{equation*}
    f_1(x,\tau) \equiv 0 \qquad \textrm{for all} \quad x < m_1(0) \frac{e^{\alpha(\chi+1)\tau} - 1}{e^{\alpha\chi \tau}},
\end{equation*}
and similarly, when $f_2^\init(y) \equiv 0$, for any $y < 0$ the distribution $f_2(y,\tau)$ will be such that
\begin{equation*}
    f_2(y,\tau) \equiv 0 \qquad \textrm{for all} \quad y < m_2(0) \frac{e^{\nu(\theta+1) \tau} - 1}{e^{\nu\theta \tau}},
\end{equation*}
implying that no populations with sizes below the obtained thresholds can emerge at any given time $\tau \geq 0$. 

%%%%%%%%%%%%%%%%%%%%%%%%%%%%%%%%%%%%%%%%%%%%%%%%%%%%%%%%%%%%%%%%%%%
\subsection{Scaling of the full system}

Next, we focus on the derivation of a Fokker--Planck model based on the interaction dynamics \eqref{eq:microscopic Boltzmann} leading to the Boltzmann system \eqref{eq:Boltzmann weak}. Introducing once again a small parameter $\eps \ll 1$ and a corresponding timescale $\tau = \eps t$, we consider a quasi-invariant regime of the interactions \eqref{eq:microscopic Boltzmann}, corresponding to very small variations that are obtained by properly rescaling all relevant parameters in terms of $\eps$. More precisely, the parameters $\beta$ and $\gamma$ characterizing the deterministic variations \eqref{eq:deterministic effects} are rescaled as
\begin{equation*}
\beta \to \eps \beta, \qquad \gamma \to \eps \gamma,
\end{equation*}
while the diffusion constants $\sigma_1, \sigma_2$ and the parameter $s_0$ appearing in the random effects \eqref{eq:random effects} are rescaled as
\begin{equation*}
\sigma_1 \to \eps \sigma_1, \qquad \sigma_2 \to \eps \sigma_2, \qquad s_0 \to \eps.
\end{equation*}
Consequently, from the microscopic interactions \eqref{eq:microscopic Boltzmann} one computes the following Taylor expansions
\begin{equation*} %\label{eq:Taylor}
    \begin{split}
        \left\langle \varphi(x')-\varphi(x) \right\rangle = (x'-x) \varphi'(x) + \frac{1}{2}(x'-x)^2 \varphi''(x) + \frac{1}{6}(x'-x)^3 \varphi'''(\hat{x}), \\[2mm]
        \left\langle \psi(y')-\psi(y) \right\rangle = (y'-y) \psi'(y) + \frac{1}{2}(y'-y)^2 \psi''(y) + \frac{1}{6}(y'-y)^3 \psi'''(\hat{y}),
    \end{split}
\end{equation*}
where $\hat{x} \in (\min\{x,x'\}, \max\{x,x'\})$ and $\hat{y} \in (\min\{y,y'\}, \max\{y,y'\})$. Plugging the above expansions
 into \eqref{eq:Boltzmann weak}, we then get
\begin{equation} \label{eq:scaled FP1}
    \begin{split}
        &\frac{\dd}{\dd \tau} \int_{\R_+} \varphi(x) f_1^\eps(x,\tau) \dd x = \frac{1}{\eps} \int_{\R_+ \times \R_+} \kappa_1(x,y) \big(-\Phi(y) x \varphi'(x) \big) f_1^\eps(x,\tau) f_2^\eps(y,\tau) \dd x \dd y \\[2mm]
        & \hspace*{3cm} + \frac{1}{2\eps}\int_{\R_+ \times \R_+} \kappa_1(x,y) \left[ \left( \eps^2\beta^2 \Phi^2(y)x^2 + \eps x^{2p} \frac{\sigma_1 y}{1+y} \right) \varphi''(x) \right] f_1^\eps(x,\tau) f_2^\eps(y,\tau) \dd x \dd y \\[6mm]
        & \hspace*{3cm} + \mathcal{R}_\varphi(f_1^\eps,f_2^\eps) + \tilde{\mathcal{R}}_\varphi(f_1^\eps,f_2^\eps), 
    \end{split}
\end{equation}
and 
\begin{equation} \label{eq:scaled FP2}
    \begin{split}
        & \frac{\dd}{\dd \tau} \int_{\R_+} \psi(y) f_2^\eps(y,\tau) \dd y = \frac{1}{\eps} \int_{\R_+ \times \R_+} \kappa_2(x,y) \big( \Psi(x) y \psi'(y) \big) f_1^\eps(x,\tau) f_2^\eps(y,\tau) \dd x \dd y \\[2mm]
        & \hspace*{3cm} + \frac{1}{2\eps} \int_{\R_+ \times \R_+} \kappa_2(x,y) \left[ \left( \eps^2\gamma^2 \Psi^2(x) y^2 + \eps y^{2p} \frac{\sigma_2 x}{1+x} \right) \psi''(y) \right] f_1^\eps(x,\tau)f_2^\eps(y,\tau) \dd x \dd y \\[6mm]
        & \hspace*{3cm} + \mathcal{R}_\psi(f_1^\eps,f_2^\eps) + \tilde{\mathcal{R}}_\psi(f_1^\eps,f_2^\eps).
    \end{split}
\end{equation}
In \eqref{eq:scaled FP1}  and \eqref{eq:scaled FP2} the remainder terms $\mathcal{R}_\varphi(f_1^\eps,f_2^\eps)$ and $\mathcal{R}_\psi(f_1^\eps,f_2^\eps)$ read
\begin{equation*}
    \begin{split}
        & \mathcal{R}_\varphi(f_1^\eps,f_2^\eps)(x,\tau) = \frac{1}{6\eps} \int_{\R_+ \times \R_+}\kappa_1(x,y) \big\langle(-\Phi(y)x + x^p \eta_1(y))^3 \big\rangle \varphi'''(\hat{x}) f_1^\eps(x,\tau) f_2^\eps(y,\tau) \dd x \dd y, \\[4mm]
        & \mathcal{R}_\psi(f_1^\eps,f_2^\eps)(y,\tau) = \frac{1}{6\eps} \int_{\R_+ \times \R_+} \kappa_2(x,y)\big\langle(\Psi(x)y + y^p \eta_2(x))^3 \big\rangle \psi'''(\hat{y}) f_1^\eps(x,\tau) f_2^\eps(y,\tau) \dd x \dd y, 
    \end{split}
\end{equation*}
while the other remainders $\tilde{\mathcal{R}}_\varphi(f_1^\eps,f_2^\eps)$ and $\tilde{\mathcal{R}}_\psi(f_1^\eps,f_2^\eps)$ are defined as
\begin{equation*}
    \begin{split}
        \tilde{\mathcal{R}}_\varphi(f_1^\eps,f_2^\eps) =& \int_{\R_+ \times \R_+}\kappa_1(x,y) \mathds{1}(0\leq x \leq (1-p)\eps) x^{2p} \left\langle \eta_1^2(y) \right\rangle \varphi''(\hat{x}) f_1^\eps(x,t) f_2^\eps(y,t) dx dy \\[2mm]
        & + \int_{\R_+ \times \R_+}\kappa_1(x,y) \mathds{1}(0 \leq x \leq (1-p)\eps) x^{3p} \left\langle 3\eta_1^2(y) + \eta_1^3(y) \right\rangle \varphi'''(\hat{x}) f_1^\eps(x,t) f_2^\eps(y,t) dx dy \\[4mm]
        \tilde{\mathcal{R}}_\psi(f_1^\eps,f_2^\eps) =& \int_{\R_+ \times \R_+}\kappa_2(x,y) \mathds{1}(0 \leq y \leq (1-p)\eps) y^{2p} \left\langle \eta_2^2(x) \right\rangle \psi''(\hat{y}) f_1^\eps(x,t) f_2^\eps(y,t) dx dy \\[2mm]
        & + \int_{\R_+ \times \R_+}\kappa_2(x,y) \mathds{1}(0 \leq y \leq (1-p)\eps) y^{3p} \left\langle 3\eta_2^2(x) + \eta_2^3(x) \right\rangle \psi'''(\hat{y}) f_1^\eps(x,t) f_2^\eps(y,t) dx dy,
    \end{split}
\end{equation*}
and take into account the presence of the indicator functions in the random effects of the elementary interactions \eqref{eq:microscopic Boltzmann}. These terms disappear in the case $p = 1$, but whenever $p < 1$ they are different from zero on a subset of measure $(1-p)\eps \ll 1$ with respect to $(x,y)$.

Now, assuming that the random variables $\eta_1(y)$ and $\eta_2(x)$ have bounded third order moments, i.e. $\left\langle |\eta_1(y)|^3 \right\rangle < +\infty$ and $\left\langle |\eta_2(x)|^3 \right\rangle <+\infty$, then we can write $\eta_1 = \sqrt{\dfrac{\sigma_1 y}{1+y}}\hat{\eta}_1$ and $\eta_2 = \sqrt{\dfrac{\sigma_2 x}{1+x}}\hat{\eta}_2$, with $\left\langle\hat{\eta}_1\right\rangle = \left\langle\hat{\eta}_2\right\rangle = 0$ and $\left\langle\hat{\eta}_1^2\right\rangle = \left\langle\hat{\eta}_2^2\right\rangle = 1$. Therefore, observing that 
\begin{equation*}
    \left| \mathcal R_\varphi(f_1^\eps,f_2^\eps)\right| \approx  \eps + \eps^2 + \sqrt{\eps}, \qquad 
    \left| \mathcal R_\psi(f_1^\eps,f_2^\eps)\right| \approx  \eps + \eps^2 + \sqrt{\eps},
\end{equation*}
and also
\begin{equation*}
    |\tilde{\mathcal{R}}_\varphi(f_1^\eps,f_2^\eps)| \approx \eps^{1+2p}, \qquad |\tilde{\mathcal R}_\psi(f_1^\eps,f_2^\eps)| \approx \eps^{1+2p},
\end{equation*}
in the limit $\eps \to 0$ all remainder terms vanish and the distributions $f_1^\eps(x,\tau)$ and $f_2^\eps(y,\tau)$ formally converge to some $f_1(x,\tau)$ and $f_2(y,\tau)$ which are solutions to the limit system of \eqref{eq:scaled FP1}--\eqref{eq:scaled FP2}, being
\begin{equation*}
    \begin{split}
        \frac{\dd}{\dd \tau} \int_{\R_+} \varphi(x) f_1(x,\tau) \dd x &= \int_{\R_+ \times \R_+} \beta x y \varphi'(x)  f_1(x,\tau) f_2(y,\tau) \dd x \dd y \\[2mm]
        & \hspace*{1cm} + \frac{\sigma_1}{2} \int_{\R_+ \times \R_+} x^{2p} y \varphi''(x) f_1(x,\tau) f_2(y,\tau) \dd x \dd y, \\[4mm]
        \frac{\dd}{\dd \tau} \int_{\R_+} \psi(y) f_2(y,\tau) \dd y &= \int_{\R_+ \times \R_+} \gamma (x - \mu) y \psi(y)' f_1(x,\tau) f_2(y,\tau) \dd x \dd y \\[2mm]
        & \hspace*{1cm} + \frac{\sigma_2}{2} \int_{\R_+ \times \R_+} y^{2p} x \psi''(y) f_1(x,\tau) f_2(y,\tau) \dd x \dd y. 
    \end{split}
\end{equation*}
Integrating back by parts, the previous formulation writes in strong form as 
\begin{equation} \label{eq:FP Boltzmann}
    \begin{split}
        \frac{\partial f_1(x,\tau)}{\partial \tau} &= \frac{\sigma_1 m_2(\tau)}{2} \frac{\partial^2}{\partial x^2} \left[ x^{2p}f_1(x,\tau) \right] + \frac{\partial}{\partial x} \big(\beta m_2(\tau) x f_1(x,\tau) \big), \\[2mm]
        \frac{\partial f_2(y,\tau)}{\partial \tau} &= \frac{\sigma_2 m_1(\tau)}{2} \frac{\partial^2}{\partial y^2} \left[ y^{2p} f_2(y,\tau) \right] + 
        \frac{\partial}{\partial y} \left[\gamma(\mu - m_1(\tau)) y f_2(y,\tau) \right] 
    \end{split}
\end{equation}
Finally, restoring the time variable $t \geq 0$ and merging the contributions \eqref{eq:FP redistribution} and \eqref{eq:FP Boltzmann} coming from the redistribution and the Boltzmann operators, we end up with the following system of Fokker--Planck-type equations
\begin{equation} \label{eq:Fokker-Planck}
    \begin{split}
        \frac{\partial f_1(x,t)}{\partial t} &= \frac{\sigma_1 m_2(t)}{2} \frac{\partial^2}{\partial x^2} \left[ x^{2p}f_1(x,t) \right] + \frac{\partial}{\partial x} \big( ((\beta m_2(t) + \alpha \chi) x -\alpha(\chi+1) m_1(t)) f_1(x,t) \big), \\[6mm]
        \frac{\partial f_2(y,t)}{\partial t} &= \frac{\sigma_2m_1(t)}{2} \frac{\partial^2}{\partial y^2} \left[y^{2p}f_2(y,t)\right] + \frac{\partial}{\partial y} \big( (\gamma(\mu - m_1(t)) y + \nu\theta y -\nu(\theta+1) m_2(t)) f_2(y,t) \big),
    \end{split}
\end{equation}
for any $0 < p \leq 1$, complemented by the no-flux boundary conditions
\begin{equation*}
    \begin{split}
        & \frac{\sigma_1 m_2(t)}{2} \frac{\partial}{\partial x} \left[x^{2p}f_1(x,t)\right] + ((\beta m_2(t) + \alpha\chi) x - \alpha(\chi+1) m_1(t)) f_1(x,t) \Big|_{x=0, +\infty} = 0, \\[2mm]
        & x^{2p} f_1(x,t) \Big|_{x=0, +\infty} = 0, \\[6mm]
        & \frac{\sigma_2 m_1(t)}{2} \frac{\partial}{\partial y} \left[y^{2p}f_2(y,t)\right] + (\gamma(\mu - m_1(t)) y + \nu\theta y - \nu(\theta+1) m_2(t)) f_2(y,t) \Big|_{y=0, +\infty} = 0, \\[2mm]
        & y^{2p} f_2(y,t) \Big|_{y=0, +\infty} = 0. 
\end{split}
\end{equation*}
Well-posedness of the above system can be assessed by applying the results from the work of Le Bris and Lions \cite{Bris03072008}.

%%%%%%%%%%%%%%%%%%%%%%%%%%%%%%%%%%%%%%%%%%%%%%%%%%%%%%%%%%%%%%%%%%%
\subsection{Evolution of the observable quantities}

It is immediate to show that the evolution of the first order moments of the solutions to system \eqref{eq:Fokker-Planck} is consistent with the one obtained from the Boltzmann-type equations \eqref{eq:Boltzmann weak}, so that the means $m_1(t)$ and $m_2(t)$ of the two populations are still solutions to the Lotka--Volterra system \eqref{eq:Lotka-Volterra}. The added value in working with this simplified Fokker--Planck formulation lies in the fact that we can now recover a closed evolution of the higher moments of the solutions to \eqref{eq:Fokker-Planck}.
%However, it is the variances, $v_1(t)$ and $v_2(t)$, that can be obtained in closed form by the Fokker--Planck system. Indeed,
%\begin{align}
%    \frac{\dd v_1(t)}{\dd t} & = -2 (\beta m_2(t) + \alpha \chi) v_1(t) + \sigma_1 m_1(t) m_2(t), \label{eq:variance-f} \\[2mm]
%    \frac{\dd v_2(t)}{\dd t} & = \textcolor{red}{-2\gamma (\mu - m_1(t))} v_2(t) + \sigma_2 m_1(t) m_2(t). \label{eq:variance-g}
%\end{align}
%These equations coincide with equations \eqref{varf}-\eqref{varg} in the limit $\eps \to 0$. This is a distinctive feature of the Fokker--Planck picture. Indeed, the system of variances \eqref{varf}-\eqref{varg}, gained at Boltzmann equations level, does not provide a closed form, as some terms, not related to any moment of the system, appear.
In particular, by concentrating on the two prototypical cases $p = 1$ and $p = 1/2$, we can determine the explicit dynamics of the variances $v_1(t)$ and $v_2(t)$. For $p = 1$ we obtain
%\begin{equation}
%\begin{split}
%    \frac{\dd E_1}{\dd t} & = -2(\beta m_2 + \alpha\chi) E_1 + 2\alpha(\chi + 1) m_1^2 + \sigma_1 m_1 m_2, \\
%    \frac{\dd E_2}{\dd t} & = -2[ \gamma(\mu - m_1) + \nu\theta] E_2 + 2\nu(\theta+1) m_2^2 + \sigma_2m_1 m_2,
%\end{split}
%\end{equation}
\begin{equation} \label{eq:variance FP p=1}
    \begin{split}
        & \frac{\dd v_1(t)}{\dd t} = - 2 \big( \big(\beta - \sigma_1/2\big) m_2(t) + \alpha \chi \big) v_1(t) + \sigma_1 m_1^2(t) m_2(t), \\[2mm]
        & \frac{\dd v_2(t)}{\dd t} = - 2 \big( \gamma(\mu - (1 - \sigma_2/2\gamma) m_1(t)) + \nu \theta \big) v_2(t) + \sigma_2 m_1(t) m_2^2(t),
    \end{split}
\end{equation}
while choosing $p = 1/2$ leads to
%\begin{equation}
%\begin{split}
%    \frac{\dd E_1(t)}{\dd t} & = -\left[ \big(2\beta-\sigma_1\big) m_2 + 2\alpha\chi \right] E_1 + 2\alpha(\chi+1) m_1^2, \\
%    \frac{\dd E_2(t)}{\dd t} & = -\left[ 2\gamma(\mu-m_1) -\sigma_2 m_1 + 2\nu\theta \right] E_2 + 2\nu(\theta+1) m_2^2,
%\end{split}
%\end{equation}
\begin{equation} \label{eq:variance FP p=1/2}
    \begin{split}
        & \frac{\dd v_1(t)}{\dd t} = -2(\beta m_2(t) + \alpha\chi) v_1(t) + \sigma_1 m_1(t) m_2(t), \\[2mm]
        & \frac{\dd v_2(t)}{\dd t} = -2 \big( \gamma(\mu - m_1(t)) + \nu\theta \big) v_2(t) + \sigma_2 m_1(t) m_2(t).
    \end{split}
\end{equation}
We highlight that these two systems governing the evolution of the variances for the Fokker--Planck equations \eqref{eq:Fokker-Planck} are in fact coherent with the ones determined from the original Boltzmann equations \eqref{eq:Boltzmann weak}, in the quasi-invariant regime of parameters considered. Furthermore, we notice that system \eqref{eq:variance FP p=1} loses stability compared to system \eqref{eq:variance FP p=1/2}, due to the diffusion parameter $\sigma_1 > 0$ reducing the region of parameters in which the variances remain bounded. A similar effect will be also observed in the following part, where we shall see that the parameter $p$ determines the structure of the equilibrium states, producing fat-tails and thin-tails distributions respectively in the cases $p = 1$ and $p = 1/2$. This loss of higher order moments when $p = 1$ will eventually lead us to consider only the more stable case $p = 1/2$.

%%%%%%%%%%%%%%%%%%%%%%%%%%%%%%%%%%%%%%%%%%%%%%%%%%%%%%%%%%%%%%%%%%%
\subsection{Study of the local equilibria}
A further property of the Fokker--Planck description is related to the possibility of obtaining the explicit form of the equilibria, which are determined by setting the left-hand sides of equations \eqref{eq:Fokker-Planck} equal to zero. In the case of the classical Fokker--Planck equations for gases \cite{MR3597010,ParTos-2013}, the knowledge of the uniform-in-time steady state and of the (exponential) rate of relaxation of the solution toward it allows to conclude that the distribution solving the Fokker--Planck system under study is well-represented by the characteristics of its corresponding equilibrium. However, the situation is quite different here, because system \eqref{eq:Fokker-Planck} does not possess stable uniform steady states, but only time-dependent oscillating local equilibria, and it is not clear whether the characteristics of these local equilibria do represent or not a good approximation of the real solutions to the system, even for intermediate and for large times. This is a challenging question which we try to partially answer here, but that it could deserve further investigations. We mention the works \cite{ARSW,Matthes} for related settings in cross-diffusion systems.

The local equilibrium states (or \emph{quasi-equilibria}) of the Fokker--Planck system \eqref{eq:Fokker-Planck} are all vector distribution functions $\mathbf{f}^\eq(\x,t) = (f_1^\eq(x,t),f_2^\eq(y,t))$ solving
\begin{equation} \label{eq:equilibrium relations}
    \begin{split}
        & \frac{\partial}{\partial x} \big(x^{2p} f_1^\eq(x,t)\big) + \frac{2}{\sigma_1 m_2(t)} \big(((\beta m_2(t) + \alpha \chi)x - \alpha (\chi+1)m_1(t)) f_1^\eq(x,t)\big) = 0, \\[2mm]
        & \frac{\partial}{\partial y} \big(y^{2p} f_2^\eq(y,t)\big) + \frac{2}{\sigma_2 m_1(t)} \big((\gamma(\mu - m_1(t)) + \nu \theta) y - \nu(\theta + 1) m_2(t)) f_2^\eq(y,t)\big) = 0,
    \end{split}
\end{equation}
for any $t > 0$ and $x,y \in \R_+$. The above system can be conveniently recast as
\begin{align*}
    & \frac{\partial}{\partial x} \log f_1^\eq(x,t) = \frac{-\frac{2}{\sigma_1 m_2(t)}((\beta m_2(t) + \alpha \chi) x - \alpha (\chi+1)m_1(t)) - 2p x^{2p-1}}{x^{2p}}, \\[4mm]
    & \frac{\partial}{\partial y} \log f_2^\eq(y,t) = \frac{-\frac{2}{\sigma_2 m_1(t)}(\gamma(\mu - m_1(t)) y - \nu m_2(t)) - 2p y^{2p-1}}{y^{2p}},
\end{align*}
and, reducing once again to the study of the main representative cases $p = 1$ and $p = 1/2$, by simple computations we find the following explicit form of the equilibrium states, for $p = 1$
\begin{align*}
    f_1^\eq(x,t) & = C_1(t) x^{- \left(2 \frac{(\beta + \sigma_1/2) m_2(t) + \alpha \chi}{\sigma_1 m_2(t)} + 1\right)} \exp\left\{-\frac{2 \alpha (\chi+1) m_1(t)}{\sigma_1 m_2(t)} \frac{1}{x} \right\}, \\[4mm]
    f_2^\eq(y,t) & = C_2(t) y^{- \left(2 \frac{\gamma(\mu - (1 - \sigma_2/2\gamma) m_1(t) + \nu \theta}{\sigma_2 m_1(t)} + 1\right)} \exp\left\{-\frac{2 \nu (\theta+1) m_2(t)}{\sigma_2 m_1(t)} \frac{1}{y} \right\},
\end{align*}
and for $p = 1/2$
\begin{equation} \label{eq:equilibrium states}
    \begin{split}
        f_1^\eq(x,t) & = C_1(t) x^{\frac{2\alpha(\chi+1)}{\sigma_1}\frac{m_1(t)}{m_2(t)} - 1} \exp\left\{-\frac{2}{\sigma_1 m_2(t)} (\beta m_2(t) + \alpha \chi) x \right\},
        \\[2mm]
        f_2^\eq(y,t) & = C_2(t) y^{\frac{2 \nu(\theta+1)}{\sigma_2}\frac{m_2(t)}{m_1(t)} - 1} \exp\left\{-\frac{2}{\sigma_2 m_1(t)} (\gamma(\mu - m_1(t) + \nu\theta) y \right\},
    \end{split}
\end{equation}
where $C_1(t)$, $C_2(t) > 0$ are normalization coefficients depending on the means $m_1(t)$ and $m_2(t)$. In particular, one recognizes the structure of the inverse Gamma distribution emerging for $p = 1$ and that of the Gamma distribution for $p = 1/2$, provided of course that we appropriately bound the parameters of the problem in order to ensure that the exponents remain positive. However, since an inverse Gamma distribution does not possess finite moments of any order, from now on we shall focus our attention on the more interesting case $p = 1/2$, corresponding to thin-tails local equilibrium states.

Let us then rewrite, for $p = 1/2$, the equilibrium distributions \eqref{eq:equilibrium states} as
\begin{equation*} %\label{eq:qe_malthus}
    f_1^\eq(x,t) = C_1(t) x^{a_1(t) - 1} e^{- b_1(t) x}, \qquad f_2^\eq(y,t) = C_2(t) y^{a_2(t) - 1} e^{- b_2(t) y},
\end{equation*}
where we have denoted with
\begin{equation*} %\label{parameters}
    \begin{split}
        & a_1(t) = \frac{2\alpha(\chi+1)}{\sigma_1}\frac{m_1(t)}{m_2(t)}, \qquad b_1(t) = \frac{2}{\sigma_1 m_2(t)} (\beta m_2(t) + \alpha \chi), \\[2mm]
        & a_2(t) = \frac{2\nu(\theta+1)}{\sigma_2}\frac{m_2(t)}{m_1(t)}, \qquad b_2(t) = \frac{2}{\sigma_2 m_1(t)} (\gamma(\mu - m_1(t)) + \nu\theta),
    \end{split}
\end{equation*}
the four exponents of interest. To recover a Gamma distribution, we must then guarantee that the latter are positive. Assuming that we have fixed the main parameters $\alpha$, $\beta$, $\delta$, $\gamma > 0$ of the Lotka--Volterra system \eqref{eq:Lotka-Volterra}, together with the interaction coefficient $\nu > 0$ and the lowest group of preys $\mu \geq 1$ which define $\delta = \gamma \mu - \nu$, the positivity of the above exponents depends on the bounds \eqref{eq:bounds on means} for the means $m_1(t)$ and $m_2(t)$. In particular, one needs to further impose that $\chi, \theta \in \R$ satisfy the inequalities
\begin{equation} \label{eq:admissible parameters}
    \chi > \max\left\{ -1,\ -\frac{\beta}{\alpha}\underline{c}_2 \right\}, \qquad \theta > \max\left\{ -1,\ -\frac{\gamma}{\nu} (\mu - \overline{c}_1) \right\},
    %\begin{array}{llll}
    %    \displaystyle \chi > - \max\left\{ 1,\ \frac{\beta}{\alpha}\underline{c}_2 \right\} & \textrm{and} & \displaystyle \theta > 0 & (\textrm{if} \ \mu \geq \overline{c}_1), \\[6mm]
    %    \displaystyle \chi > - \max\left\{ 1,\ \frac{\beta}{\alpha}\underline{c}_2 \right\} & \textrm{and} & \displaystyle \theta > \frac{\gamma}{\nu} (\overline{c}_1 - \mu) & (\textrm{if} \ \underline{c}_1 \leq \mu < \overline{c}_1).
    %\end{array}
\end{equation}
Under these bounds, we can exploit the properties of the Gamma distributions to determine an explicit expression for the normalization constants $C_1(t)$ and $C_2(t)$, which read
\begin{equation*}
    C_1(t) = \frac{b_1^{a_1}}{\Gamma(a_1)}, \qquad C_2(t) = \frac{b_2^{a_2}}{\Gamma(a_2)}.
\end{equation*}
Moreover, by denoting with $\zeta_1 > 0$ and $\zeta_2 > 0$ the two positive lower bounds for which conditions \eqref{eq:admissible parameters} are satisfied, namely $\beta \underline{c}_2 + \alpha \chi \geq \zeta_1$ and $\gamma(\mu - \overline{c}_1) + \nu\theta \geq \zeta_2$, one can compute the explicit solutions of system \eqref{eq:variance FP p=1/2} as
\begin{equation} \label{eq:variances explicit}
    \begin{split}
        v_1(t) &= v_1(0) e^{-2 \int_0^t (\beta m_2(s) + \alpha \chi) \dd s} + \sigma_1 \int_0^t m_1(s) m_2(s) e^{-2 \int_s^t (\beta m_2(\bar{s}) + \alpha \chi) \dd \bar{s}} \dd s, \\[2mm]
        v_2(t) &= v_2(0) e^{-2 \int_0^t (\gamma(\mu - m_1(s)) + \nu\theta) \dd s} + \sigma_2 \int_0^t m_1(s) m_2(s) e^{-2 \int_s^t (\gamma(\mu - m_1(\bar{s})) + \nu\theta) \dd \bar{s}} \dd s,
    \end{split}
\end{equation} 
and consequently infer the following global-in-time $L^\infty$ controls for the variances $v_1(t)$ and $v_2(t)$
\begin{equation} \label{eq:bounds on variances}
    \begin{split}
        \| v_1 \|_{L^\infty(\R_+)} &\leq v_1(0) e^{-2 \zeta_1 t} + \frac{\sigma_1}{2 \zeta_1} \overline{c}_1 \overline{c}_2 \left( 1 - e^{-2 \zeta_1 t}\right), \\[2mm]
        \| v_2 \|_{L^\infty(\R_+)} &\leq v_2(0) e^{-2 \zeta_2 t} + \frac{\sigma_2}{2 \zeta_2} \overline{c}_1 \overline{c}_2 \left( 1 - e^{-2 \zeta_2 t}\right),
    \end{split}
\end{equation}
based on the global bounds $\| m_1 \|_{L^\infty(\R_+)} = \overline{c}_1$ and $\| m_2 \|_{L^\infty(\R_+)} = \overline{c}_2$ for the means.

\begin{remark}
We highlight that sharper a priori estimates show that the variances are actually well-defined for ranges of admissible $\chi$ and $\theta$ that are wider than the ones imposed by conditions \eqref{eq:admissible parameters}. Indeed, exploiting the integral form of the solutions to the Lotka--Volterra system \eqref{eq:Lotka-Volterra}, namely
\begin{equation*}
    \begin{split}
        m_1(t) &= m_1(0) e^{-\int_0^t (\beta m_2(s) - \alpha) \dd s}, \\[4mm]
        m_2(t) &= m_2(0) e^{-\int_0^t (\delta - \gamma m_1(s)) \dd s},
    \end{split}
\end{equation*}
one can extract these terms from system \eqref{eq:variances explicit} by adding and subtracting respectively $\alpha$ and $\nu$ inside the integrals that define the exponents on the right-hand sides. In doing so, simple algebraic manipulations allow to recast \eqref{eq:variances explicit} as
\begin{equation*}
    \begin{split}
        v_1(t) &= \frac{v_1(0)}{m_1^2(0)} m_1^2(t) e^{-2 \alpha(\chi+1) t} + \sigma_1 m_1^2(t) e^{-2 \alpha(\chi+1) t} \int_0^t \frac{m_2(s)}{m_1(s)} e^{2 \alpha(\chi+1) s} \dd s, \\[2mm]
        v_2(t) &= \frac{v_2(0)}{m_2^2(0)} m_2^2(t) e^{-2  \nu(\theta+1) t} + \sigma_2 m_2^2(t) e^{-2 \nu(\theta+1) t} \int_0^t \frac{m_1(s)}{m_2(s)} e^{-2 \nu(\theta+1) s} \dd s,
    \end{split}
\end{equation*}
from which we derive the global-in-time controls
\begin{equation*}
    \begin{split}
        \| v_1 \|_{L^\infty(\R_+)} &\leq \frac{v_1(0)}{m_1^2(0)} \overline{c}_1^2 e^{-2 \alpha(\chi+1) t} + \frac{\sigma_1}{2 \alpha(\chi+1)} \frac{\overline{c}_2}{\underline{c}_1} \left( 1 - e^{-2 \alpha(\chi+1) t}\right), \\[2mm]
        \| v_2 \|_{L^\infty(\R_+)} &\leq \frac{v_2(0)}{m_2^2(0)} \overline{c}_2^2 e^{-2  \nu(\theta+1) t} + \frac{\sigma_2}{2 \nu(\theta+1)} \frac{\overline{c}_1}{\underline{c}_2} \left( 1 - e^{-2 \nu(\theta+1) t}\right),
    \end{split}
\end{equation*}
implying that the variances remain in fact bounded as long as $\chi > -1$ and $\theta > -1$.
\end{remark}

%%%%%%%%%%%%%%%%%%%%%%%%%%%%%%%%%%%%%%%%%%%%%%%%%%%%%%%%%%%%%%%%%%%
\subsection{Asymptotic behavior of the moments}

In this part, we analyze how far the solution $\mathbf{f}(\x,t) = (f_1(x,t),f_2(y,t))$ to the Fokker--Planck system \eqref{eq:Fokker-Planck} is from the local Gamma equilibrium $\mathbf{f}^\eq(\x,t) = (f_1^\eq(x,t), f_2^\eq(y,t))$ defined by \eqref{eq:equilibrium states}, in terms of the distances between means and between variances of the two vector-distributions.

Given this equilibrium state $\mathbf{f}^\eq(\x,t)$, we define the vector $\mathbf{m}^\eq(t) = (m_1^\eq(t) ,m_2^\eq(t))$ of the corresponding means $\displaystyle m_1^\eq(t) = \int_{\R_+} x f_1^\eq(x,t) \dd x$ and $\displaystyle m_2^\eq(t) = \int_{\R_+} y f_2^\eq(y,t) \dd y$. The explicit form of these means can be determined thanks to well-known formulae for the Gamma distributions, specifically
\begin{equation} \label{eq:Gamma means}
    \begin{split}
        m_1^\eq(t) &= \frac{a_1(t)}{b_1(t)} = \frac{\alpha(\chi+1)m_1(t)}{\beta m_2(t) + \alpha \chi}, \\[4mm]
        m_2^\eq(t) &= \frac{a_2(t)}{b_2(t)} = \frac{\nu(\theta+1) m_2(t)}{\gamma (\mu - m_1(t)) + \nu\theta},
    \end{split}
\end{equation}
and one can measure their distance from the solution $\mathbf{m}(t) = (m_1(t), m_2(t))$ of the Lotka--Volterra system \eqref{eq:Lotka-Volterra}, by estimating the expression
\begin{equation*}
    \| \mathbf{m}(t) - \mathbf{m}^\eq(t) \|_{l^\infty} =  \frac{\beta \left| m_2(t) - \frac{\alpha}{\beta} \right|}{\beta m_2(t) + \alpha \chi}\ m_1(t) + \frac{\gamma\left| m_1(t) - \frac{\delta}{\gamma} \right|}{\gamma(\mu - m_1(t)) + \nu\theta}\ m_2(t),
\end{equation*}
where we have used that $\delta = \gamma \mu - \nu$, and the denominators are positive because of the assumptions \eqref{eq:admissible parameters} on the parameters $\chi$ and $\theta$. To recover explicit lower and upper bounds, we thus need to properly control the numerators. We first observe that since the solution $\mathbf{m}(t)$ evolves over closed orbits around the non-zero equilibrium point \eqref{eq:m*}, for any given positive initial condition $(m_1(0), m_2(0)) \neq \mathbf{m}^*$, the solution $\mathbf{m}(t)$ differs from the null vector for any $t \geq 0$, and we thus deduce the existence of a strictly positive lower bound. In particular, for any fixed $(m_1(0), m_2(0)) \in \R_+^* \times \R_+^*$, the following constants
\begin{align*}
    & \underline{c}_0 := \min_{H(m_1(t),m_2(t)) = H(m_1(0),m_2(0))} \| \mathbf{m}(t) - \mathbf{m}^* \|_{l^\infty}, \\[2mm]
    & \overline{c}_0 := \max_{H(m_1(t),m_2(t)) = H(m_1(0),m_2(0))} \| \mathbf{m}(t) - \mathbf{m}^* \|_{l^\infty},
\end{align*}
are well-defined, positive and depend only on $(m_1(0), m_2(0))$. Recalling once again the bounds \eqref{eq:bounds on means}, we may then infer the lower and upper estimates
\begin{equation} \label{eq:main estimate}
    0 < \underline{c}_0 \min \left\{ \frac{\beta \underline{c}_1}{\beta \overline{c}_2 + \alpha \chi}, \frac{\gamma \underline{c}_2}{\gamma(\mu - \underline{c}_1) + \nu\theta} \right\} \leq \| \mathbf{m}(t) - \mathbf{m}^\eq(t) \|_{l^\infty} \leq \overline{c}_0 \max \left\{ \frac{\beta \overline{c}_1}{\beta \underline{c}_2 + \alpha \chi}, \frac{\gamma \overline{c}_2}{\gamma(\mu - \overline{c}_1) + \nu\theta} \right\}.
\end{equation}
This proves existence of a positive gap in the distance over time between $\mathbf{m}(t)$ and $\mathbf{m}^\eq(t)$, preventing any possible relaxation of the solutions to the Lotka--Volterra system \eqref{eq:Lotka-Volterra} toward the means of the local equilibrium states \eqref{eq:equilibrium states} computed from the Fokker--Planck equations \eqref{eq:Fokker-Planck}, thus implying that the distributions $f_1(x,t)$ and $f_2(y,t)$ themselves cannot converge to the corresponding local equilibria $f_1^\eq(x,t)$ and $f_2^\eq(y,t)$ prescribed by the kinetic system.

%%%%%%%%%%%%%%%%%%%%%%%%%%%%%
%%%%%%%%%%%%%%%%%%%%%%%%%%%%%
\begin{table}[h!]
\begin{center}
\begin{tabular}{|c|c|c|}
\hline
Parameter		&	Value		&	Meaning \\[1mm]
\hline \hline
$\alpha$		&	1			&	Preys' growth rate \\[1mm]
$\beta$			&	0.5			&	Preys' predation rate \\[1mm]
$\mu$			& 	10			&	Preys' lowest size \\[1mm]
$\gamma$		&	0.15		&	Predators' growth rate \\[1mm]
$\nu$			&	1			&	Predators' birth rate \\[1mm]
$\delta$		& $\gamma\mu - \nu$ &	Predators' death rate \\[1mm]
$\sigma_1$		&	$10^{-3}$	&	Preys' diffusion \\[1mm]
$\sigma_2$&	$10^{-3}$	&	Predators' diffusion \\[1mm]
$\chi$			&	0			&	Preys' redistribution
                                    coefficient \\[1mm]
$\theta$		&	0			&	Predators' redistribution
                                    coefficient \\[1mm]
\hline
\end{tabular}
\end{center}
\caption{Parameters used to solve the coupled differential systems of the means \eqref{eq:Lotka-Volterra} and the variances \eqref{eq:variance FP p=1/2}.}
\label{tab:parameters}
\end{table}
%%%%%%%%%%%%%%%%%%%%%%%%%%%%%
%%%%%%%%%%%%%%%%%%%%%%%%%%%%%

In particular, we observe that $\mathbf{m}(t)$ and $\mathbf{m}^\eq(t)$ live on different closed orbits (see Figure \ref{fig:orbits of means} below) around the same equilibrium point $\mathbf{m}^*$ given by \eqref{eq:m*}, as one can easily verify by checking that $H(m_1(t),m_2(t)) \neq H(m_1^\eq(t),m_2^\eq(t))$. In Figure \ref{fig:orbits of means}, on the left, we display the shape of the orbits spanned by the solutions $\mathbf{m}(t)$ to the Lotka--Volterra system \eqref{eq:Lotka-Volterra}, for five different values of the initial condition $(m_1(0),m_2(0))$. On the right, we show the corresponding orbits of the means $\mathbf{m}^\eq(t)$ given by \eqref{eq:Gamma means}. The parameters used to run the simulations are presented in Table \ref{tab:parameters}. The above lower and upper bounds for the distance between means are therefore not surprising and one should be careful when comparing the solution to the Lotka--Volterra system \eqref{eq:Lotka-Volterra} with the means of the local Gamma equilibrium states \eqref{eq:equilibrium states}. The link between the distributions $\mathbf{f}(\x,t)$ and $\mathbf{f}^\eq(\x,t)$ can be assessed by noting that relations \eqref{eq:Gamma means} can be rewritten as
\begin{equation} \label{eq:rescaled means}
    \begin{split}
        \left( \frac{\beta m_2(t) - \alpha}{\alpha(\chi+1)} + 1 \right) m_1^\eq(t) &= m_1(t), \\[4mm]
        \left( \frac{\delta - \gamma m_1(t)}{\nu(\theta+1)} + 1 \right) m_2^\eq(t) &= m_2(t),
        \end{split}
\end{equation}
from which we deduce that rescaling the means $m_1^\eq(t)$ and $m_2^\eq(t)$ via the transformations on the left-hand sides of the above equations allows to recover the correct orbit of evolution of $\mathbf{m}(t)$.

%%%%%%%%%%%%%%%%%%%%%%%%%%%%%
%%%%%%%%%%%%%%%%%%%%%%%%%%%%%
\begin{figure}
\centering  
    \includegraphics[scale=0.4]{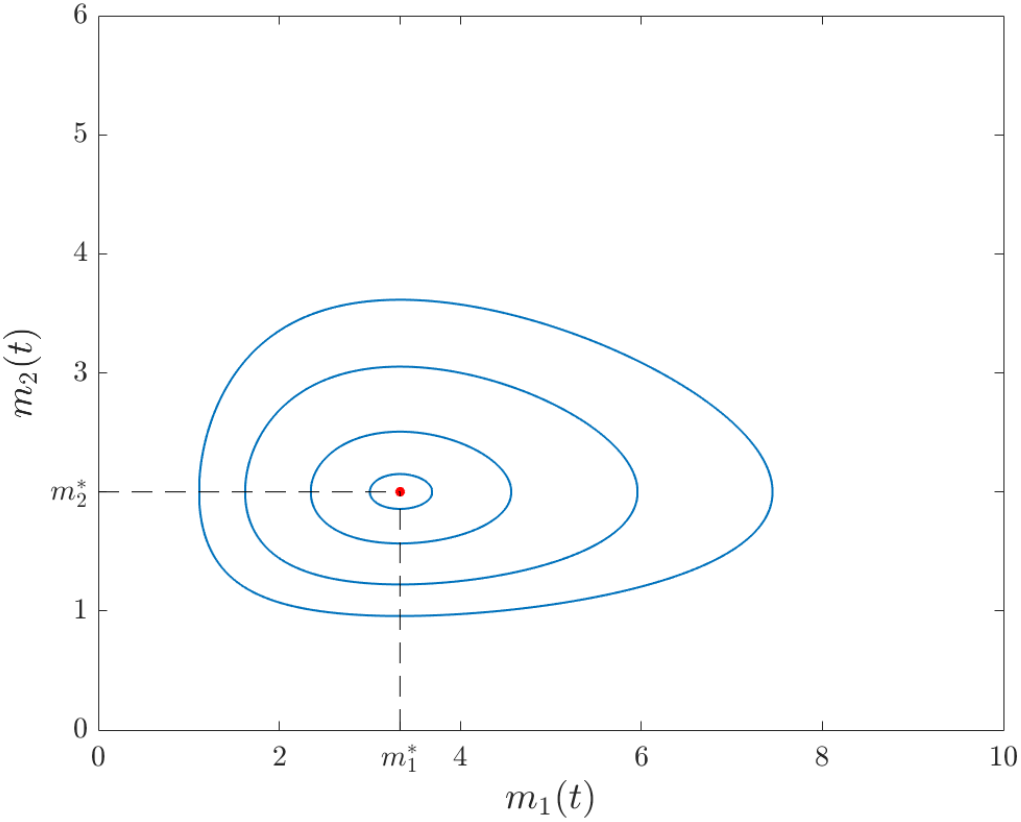}
    \hspace*{5mm}
    \includegraphics[scale=0.4]{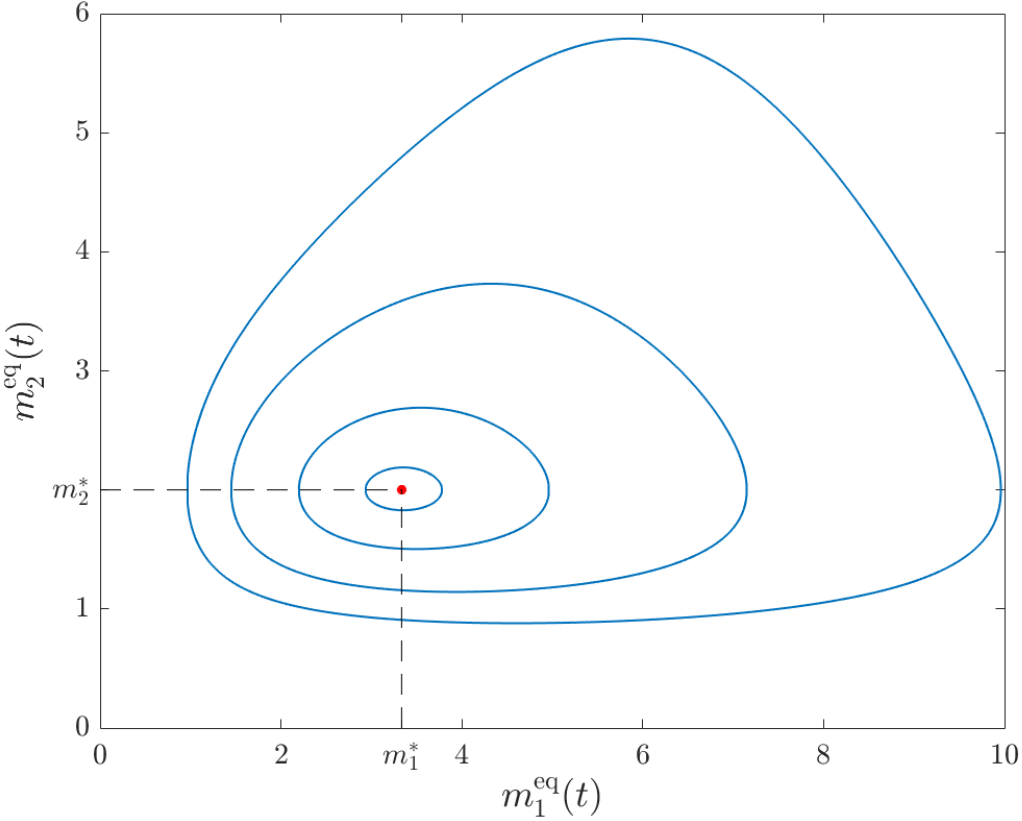}
    \caption{Orbits of evolution around the equilibrium point $\mathbf{m}^* = \left(\frac{10}{3},2\right)$ for the means $\mathbf{m}(t)$ (left), which are solutions to the Lotka--Volterra system \eqref{eq:Lotka-Volterra}, and for the means $\mathbf{m}^\eq(t)$ (right), which are computed from the local equilibrium states \eqref{eq:equilibrium states} of the Fokker--Planck equations \eqref{eq:Fokker-Planck}. The orbits (ordered by increasing size) correspond to five choices of the initial value of $\mathbf{m}(t)$, being $(m_1(0),m_2(0)) = (3,2)$, $(3.5,2.5)$, $(4,3)$, $(4.5,3.5)$. The parameters are taken from Table \ref{tab:parameters}.}
\label{fig:orbits of means}
\end{figure}
%%%%%%%%%%%%%%%%%%%%%%%%%%%%%
%%%%%%%%%%%%%%%%%%%%%%%%%%%%%

A similar behavior is also observed for the variances $\mathbf{v}(t) = (v_1(t),v_2(t))$ of the solutions $\mathbf{f}(\x,t)$ to system \eqref{eq:variance FP p=1/2}, when compared to the variances $\mathbf{v}^\eq(t) = (v_1^\eq(t),v_2^\eq(t))$ of the corresponding equilibrium state $\mathbf{f}^\eq(\x,t)$, which are defined by $\displaystyle v_1^\eq(t) = \int_{\R_+} (x - m_1^\eq(t))^2 f_1^\eq(x,t) \dd x$ and $\displaystyle v_2^\eq(t) = \int_{\R_+} (y - m_2^\eq(t))^2 f_2^\eq(y,t) \dd y$, and can be explicitly computed from the Gamma distributions \eqref{eq:equilibrium states} as
\begin{equation} \label{eq:Gamma variances}
    \begin{split}
        v_1^\eq(t) &= \frac{a_1(t)}{b_1^2(t)} = \frac{\sigma_1 \alpha(\chi+1) m_1(t) m_2(t)}{2 (\beta m_2(t) + \alpha \chi)^2}, \\[2mm]
        v_2^\eq(t) &= \frac{a_2(t)}{b_2^2(t)} = \frac{\sigma_2 \nu(\theta+1) m_1(t) m_2(t)}{2 (\gamma (\mu - m_1(t)) + \nu\theta)^2}.
    \end{split}
\end{equation}

Specifically, we notice that the right-hand sides of \eqref{eq:variance FP p=1/2} cancel out if and only if
\begin{equation*}
    \begin{split}
        v_1(t) = \frac{\sigma_1 m_1(t) m_2(t)}{2 (\beta m_2(t) + \alpha \chi)}, \qquad v_2(t) = \frac{\sigma_2 m_1(t) m_2(t)}{2 (\gamma (\mu - m_1(t)) + \nu\theta)},
    \end{split}
\end{equation*}
from which follows that the only possible equilibrium point $\mathbf{v}^*$ for the variances $\mathbf{v}(t)$ is the one determined by replacing $\mathbf{m}^*$ into the above relations. Simple computations show that
\begin{equation} \label{eq:v*}
    \mathbf{v}^* = \left(\frac{\delta \sigma_1}{2\beta \gamma(\chi +1)}, \frac{\alpha \delta \sigma_2}{2\beta \gamma \nu (\theta+1)} \right).
\end{equation}
Moreover, by rewriting the variances \eqref{eq:Gamma variances} in terms of the equilibrium means $\mathbf{m}^\eq(t)$ as
\begin{equation*}
    \begin{split}
        v_1^\eq(t) &= \frac{\sigma_1 m_1^\eq(t) m_2(t)}{2 (\beta m_2(t) + \alpha \chi)}, \\[2mm]
        v_2^\eq(t) &= \frac{\sigma_2 m_1(t) m_2^\eq(t)}{2 (\gamma (\mu - m_1(t)) + \nu\theta)},
    \end{split}
\end{equation*}
one easily deduces that $\mathbf{v}^*$ is also the unique equilibrium point of $\mathbf{v}^\eq(t)$, since $\mathbf{m}^*$ is the only equilibrium point common to both $\mathbf{m}(t)$ and $\mathbf{m}^\eq(t)$.

\begin{figure}
    \centering
    \includegraphics[scale=0.4]{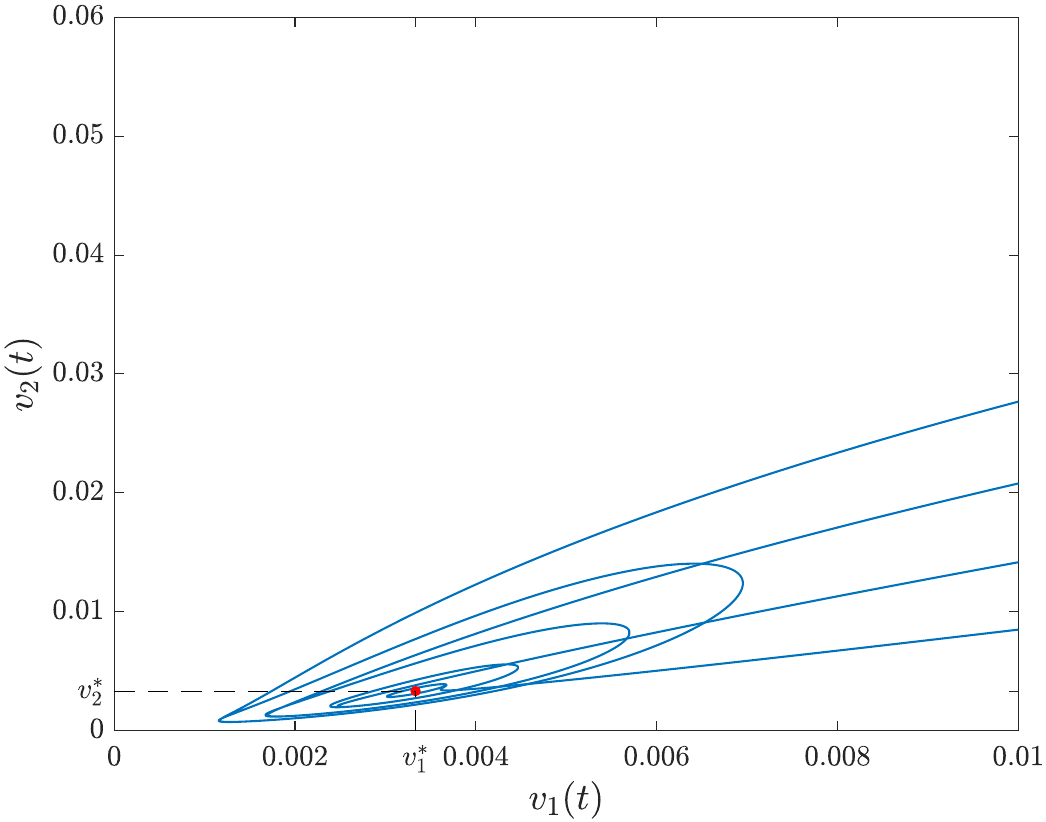}
    \hspace*{5mm}
    \includegraphics[scale=0.4]{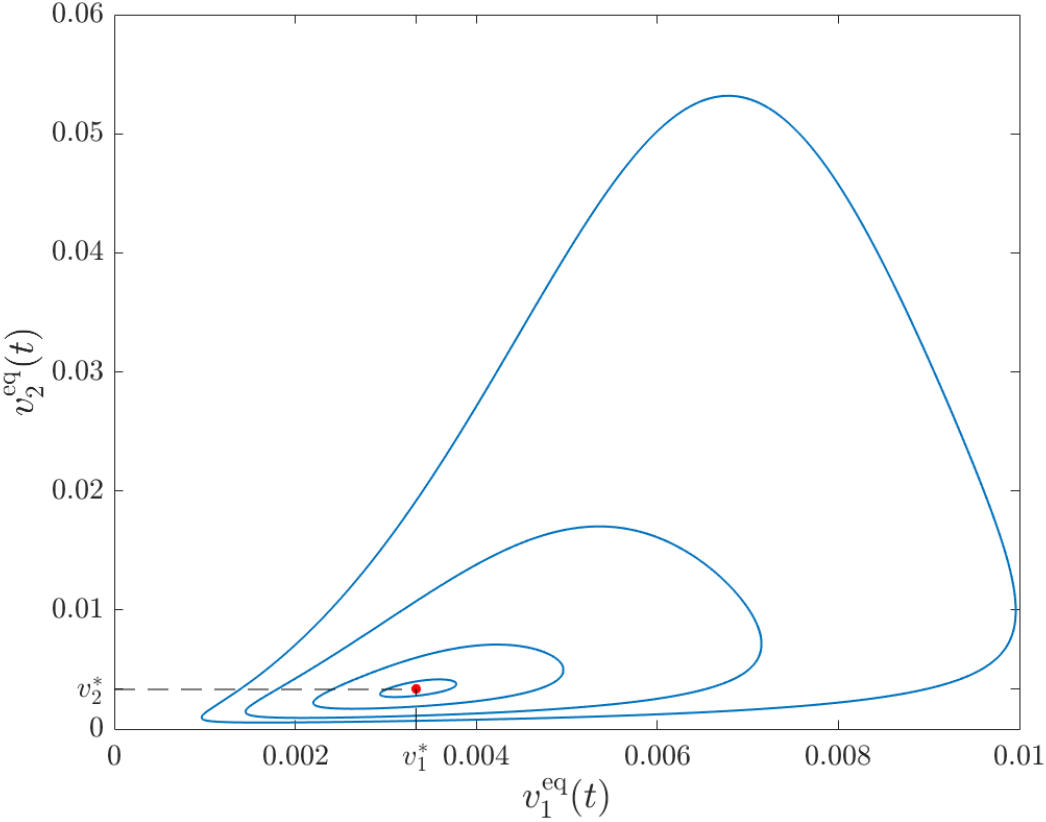}
    \caption{Left: orbits of the variances $v_1(t)$ and $v_2(t)$ defined by equations \eqref{eq:variance FP p=1/2}. Right: orbits of the quasi-equilibrium variances $v_1^\eq(t)$ and $v_2^\eq(t)$ defined by \eqref{eq:Gamma variances}. In both cases, we depict in red the point $\mathbf{v}^*$ defined in \eqref{eq:v*}, the parameters have been defined in Table \ref{tab:parameters}, the initial mean numbers of predators and preys are $(m_1(0),m_2(0)) = (3,2),(3.5,2.5),(4,3),(4.5,3.5)$ and the initial conditions for the variances are $v_1(0) = v_2(0) = 0.1$.}
    \label{fig:orbits of variances}
\end{figure}

Therefore, similar to how the means $\mathbf{m}(t)$ and $\mathbf{m}^\eq(t)$ orbit around the point $\mathbf{m}^*$, the variances $\mathbf{v}(t)$ and $\mathbf{v}^\eq(t)$ evolve, at least for large enough times, over different periodic orbits around the common equilibrium point $\mathbf{v}^*$. In fact, from the bounds \eqref{eq:bounds on variances} on the variances and by Poincaré--Bendixson theorem, one infers that the $\omega$-limit set of $\mathbf{v}(t)$ must be a periodic orbit, since the only fixed point $\mathbf{v}^*$ of system \eqref{eq:variance FP p=1/2} is never reached unless also $\mathbf{m}(t)$ equals $\mathbf{m}^*$. In particular, from the explicit expressions \eqref{eq:variances explicit} of $\mathbf{v}(t)$ we see that they decay exponentially from their initial values $v_1(0)$ and $v_2(0)$ toward a periodic orbit prescribed by the integral terms on the right-hand sides. Additionally, from the oscillating nature of the means $\mathbf{m}(t)$, one similarly conclude that the variances $\mathbf{v}^\eq(t)$ given by \eqref{eq:Gamma variances} live on periodic orbits from the start of their evolution. Unfortunately, we were not able to prove analytically the existence of a positive lower estimate similar to \eqref{eq:main estimate} for the distance $\|\mathbf{v}(t) - \mathbf{v}^\eq(t)\|_{l^\infty}$, but we can only ensure existence of a constant $C > 0$ providing an upper bound $\|\mathbf{v}(t) - \mathbf{v}^\eq(t)\|_{l^\infty} \leq C$, which follows straightforwardly from the fact that both $\mathbf{v}(t)$ and $\mathbf{v}^\eq(t)$ are bounded because of \eqref{eq:bounds on variances} and \eqref{eq:bounds on means}. We however explore these considerations numerically and depict in Figure \ref{fig:orbits of variances} an example of evolution for the orbits of the variances $\mathbf{v}(t)$ on the left and $\mathbf{v}^\eq(t)$ on the right. This confirms the lack of relaxation of $\mathbf{f}(\x,t)$ toward $\mathbf{f}^\eq(\x,t)$.

%Similarly, owing to the fact that the equilibrium state $\mathbf{f}^\eq(\x,t)$ is a vector of Gamma distributions  with parameters given by \eqref{parameters}, we can easily compute the vector of variances $\mathbf v^{\textrm{eq}}(t)  = (v_1^{\textrm{eq}}(t),v_2^{\textrm{eq}}(t))$.  We get
%\be\label{eq:Gamma variance}
%\begin{split}
%v_1^{\textrm{eq}}(t) = \dfrac{a_1(t)}{b_1^2(t)} = \dfrac{\alpha (\chi+1) \sigma_1 m_1(t) m_2(t)}{2(\beta m_2(t) + \alpha \chi)^2} \\
%v_2^{\textrm{eq}}(t) = \dfrac{a_2(t)}{b_2^2(t)} = \dfrac{\nu (\theta+1) \sigma_2 m_1(t) m_2(t)}{2(\gamma(\mu- m_2(t)) + \nu \theta)^2}. 
%\end{split}
%\ee
%Since relationships \eqref{eq:Gamma means} allow to express the denominators of relationships \eqref{eq:Gamma variance}  in terms of the quotients $m_i(t)/m_i^{eq}(t)$, $i =1,2$, we can rewrite \eqref{eq:Gamma variance} in the form
%\be\label{eq:Gamma variance2}
%\begin{split}
%v_1^{\textrm{eq}}(t) =  \dfrac{\sigma_1 [m_1^{eq}(t)]^2 m_2(t)}{2\alpha (\chi+1)m_1(t)}, \qquad v_2^{\textrm{eq}}(t) =  \dfrac{\sigma_2 [m_2^{eq}(t)]^2 m_1(t)}{2\nu (\theta+1)m_2(t)}. 
%\end{split}
%\ee
%Hence,  since both the solution $\mathbf{m}(t)$ and $\mathbf{m}^{\eq}(t)$ evolve over closed orbits around the same non-zero equilibrium point \eqref{eq:m*}. Furthermore, the evolution of the variances evolves around the same non-zero equilibrium point $\mathbf v^* $ defined in \eqref{fix-var}. In Figure \ref{fig:orbits of variances} we depict the evolution of the variances and of the variances of the quasi-equilibrium distributions.

%%%%%%%%%%%%%%%%%%%%%%%%%%%%%
%%%%%%%%%%%%%%%%%%%%%%%%%%%%%
%%%%%%%%%%%%%%%%%%%%%%%%%%%%%
%%%%%%%%%%%%%%%%%%%%%%%%%%%%%

We then draw the following probabilistic interpretation of these findings. modeling the evolution of preys and predators via the Fokker--Planck system \eqref{eq:Fokker-Planck} amounts to assume that the two populations are described by some random variables $X(t) \sim f_1(x,t)$ and $Y(t) \sim f_2(y,t)$, whose densities are given by the solution $\mathbf{f}(\x,t)$ to \eqref{eq:Fokker-Planck}. In particular, the means $m_1(t)$ and $m_2(t)$ characterize the expected values of the populations. Therefore, describing the orbits of $\mathbf{m}(t)$ through the equilibrium states \eqref{eq:equilibrium states} corresponds to looking at the evolution of the rescaled means \eqref{eq:rescaled means} and applying the same transformations to $X(t)$ and $Y(t)$, namely
\begin{equation} \label{eq:rescaled populations}
    \begin{split}
        \tilde{X}(t) = \left( \frac{\beta m_2(t) - \alpha}{\alpha(\chi+1)} + 1 \right) X(t), \\[4mm]
        \tilde{Y}(t) = \left( \frac{\delta - \gamma m_1(t)}{\nu(\theta+1)} + 1 \right) Y(t).
    \end{split}
\end{equation}
Notice in particular that the rescaling \eqref{eq:rescaled means} preserves the positivity of the means, thanks to our assumptions \eqref{eq:admissible parameters} on the admissible range of the parameters of the problem, and also the macroscopic equilibrium $\mathbf{m}^*$, since both coefficients involved in the transformation equal $1$ when $m_1(t) = \frac{\delta}{\gamma}$ and $m_2(t) = \frac{\alpha}{\beta}$. 
%when $m_1(t) = \frac{\delta}{\gamma}$ and $m_2(t) = \frac{\alpha}{\beta}$.

If we now use the Gamma distributions $\mathbf{f}^\eq(\x,t)$ to describe the evolution of $\tilde{X}(t)$ and $\tilde{Y}(t)$ at equilibrium, an easy computation shows that their expected values are given by
\begin{equation*}
    \begin{split}
        \tilde{m}_1^\eq(t) &= \int_{\R_+} \left( \frac{\beta m_2(t) - \alpha}{\alpha(\chi+1)} + 1 \right) x f_1^\eq(x,t) \dd x = \frac{\beta m_2(t) + \alpha\chi}{\alpha(\chi+1)} m_1^\eq(t),
        \\[4mm]
        \tilde{m}_2^\eq(t) &= \int_{\R_+} \left( \frac{\delta - \gamma m_1(t)}{\nu(\theta+1)} + 1 \right) y f_2^\eq(y,t) \dd y = \frac{\gamma (\mu - m_1(t)) + \nu\theta}{\nu(\theta+1)} m_2^\eq(t),
    \end{split}
\end{equation*}
which, from relations \eqref{eq:Gamma means}, equal the expected values $m_1(t)$ and $m_2(t)$ of the original variables $X(t)$ and $Y(t)$, in terms of the solution $\mathbf{f}(\x,t)$ to the Fokker--Planck system \eqref{eq:Fokker-Planck}. In this sense, the oscillating nature of the solutions to the Lotka--Volterra equations \eqref{eq:Lotka-Volterra} is correctly captured by the local equilibrium states $\mathbf{f}^\eq(\x,t)$, when the two populations are properly rescaled as in \eqref{eq:rescaled populations}.

%%%%%%%%%%%%%%%%%%%%%%%%%%%%%
%%%%%%%%%%%%%%%%%%%%%%%%%%%%%
\begin{figure}
    \centering  
    \includegraphics[scale=0.4]{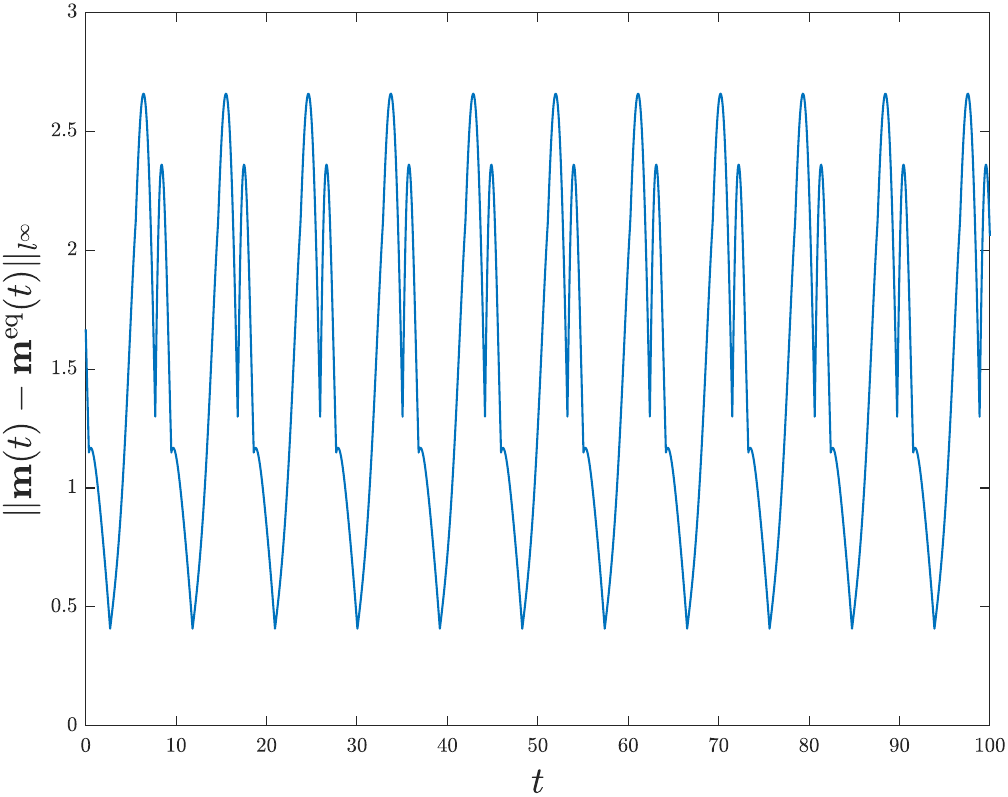}
    \hspace*{5mm}
    \includegraphics[scale=0.4]{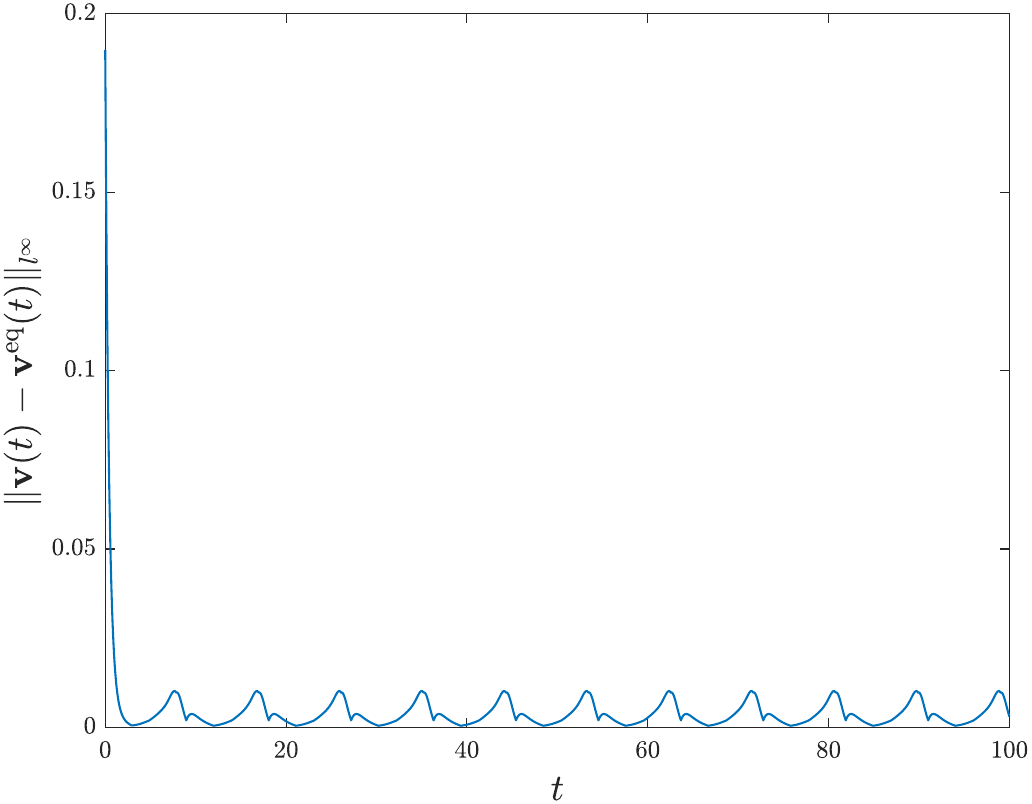}
    \caption{Analysis of the distance between the solution $\mathbf{f}(\x,t)$ to the Fokker--Planck system \eqref{eq:Fokker-Planck} and the local Gamma equilibrium state $\mathbf{f}^\eq(\x,t)$ given by \eqref{eq:equilibrium states}, in terms of their moments of order one and two. Evolution over time of the distances between the respective means (left) and variances (right). We run the simulations for a time $t \in [0,100]$, with initial conditions $(m_1(0), m_2(0)) = (4,3)$ and $(v_1(0), v_2(0)) = (0.1, 0.1)$.}
\label{fig:distance from equilibrium}
\end{figure}
%%%%%%%%%%%%%%%%%%%%%%%%%%%%%
%%%%%%%%%%%%%%%%%%%%%%%%%%%%%

In Figure \ref{fig:distance from equilibrium}, on the right, we show the evolution of the distance $\|\mathbf{v}(t) - \mathbf{v}^\eq(t)\|_{l^\infty}$ between variances, computed by solving systems \eqref{eq:Lotka-Volterra} and \eqref{eq:variance FP p=1/2} coupled together. In particular, we notice that this quantity oscillates and never reaches zero. On the left, we also plot the evolution of the distance $\|\mathbf{m}(t) - \mathbf{m}^\eq(t)\|_{l^\infty}$ between means, presenting a similar behavior that confirms our previous analytical study. The choices for the different parameters are shown in Table \ref{tab:parameters}.

Going back to the previous argument used to link the three means $\mathbf{m}(t)$, $\mathbf{m}^\eq(t)$ and $\tilde{\mathbf{m}}^\eq(t)$, we can also compare the evolution of $\mathbf{v}(t)$ with that of $\tilde{\mathbf{v}}^\eq(t)= (\tilde{v}_1^\eq(t),\tilde{v}_2^\eq(t))$, defining the variances of the rescaled populations \eqref{eq:rescaled populations} with respect to the local equilibrium $\mathbf{f}^\eq(\x,t)$, which write
\begin{equation*}
    \begin{split}
        \tilde{v}_1^\eq(t) &= \int_{\R_+} \left( \frac{\beta m_2(t) - \alpha}{\alpha(\chi+1)} + 1 \right)^2 \left(x - m_1^\eq(t)\right)^2 f_1^\eq(x,t) \dd x = \left(\frac{\beta m_2(t) + \alpha\chi}{\alpha(\chi+1)}\right)^2 v_1^\eq(t) \\[4mm]
        &= \frac{\sigma_1 m_1(t) m_2(t)}{2 \alpha(\chi+1)},
        \\[6mm]
        \tilde{v}_2^\eq(t) &= \int_{\R_+} \left( \frac{\delta - \gamma m_1(t)}{\nu(\theta+1)} + 1 \right) \left(y - m_2^\eq(t)\right)^2 f_2^\eq(y,t) \dd y = \left(\frac{\gamma (\mu - m_1(t)) + \nu\theta}{\nu(\theta+1)}\right)^2 v_2^\eq(t) \\[4mm]
        &= \frac{\sigma_2 m_1(t) m_2(t)}{2 \nu(\theta+1)}.
    \end{split}
\end{equation*}

%%%%%%%%%%%%%%%%%%%%%%%%%%%%%
%%%%%%%%%%%%%%%%%%%%%%%%%%%%%
\begin{figure}
\centering  
    \includegraphics[scale=0.3]{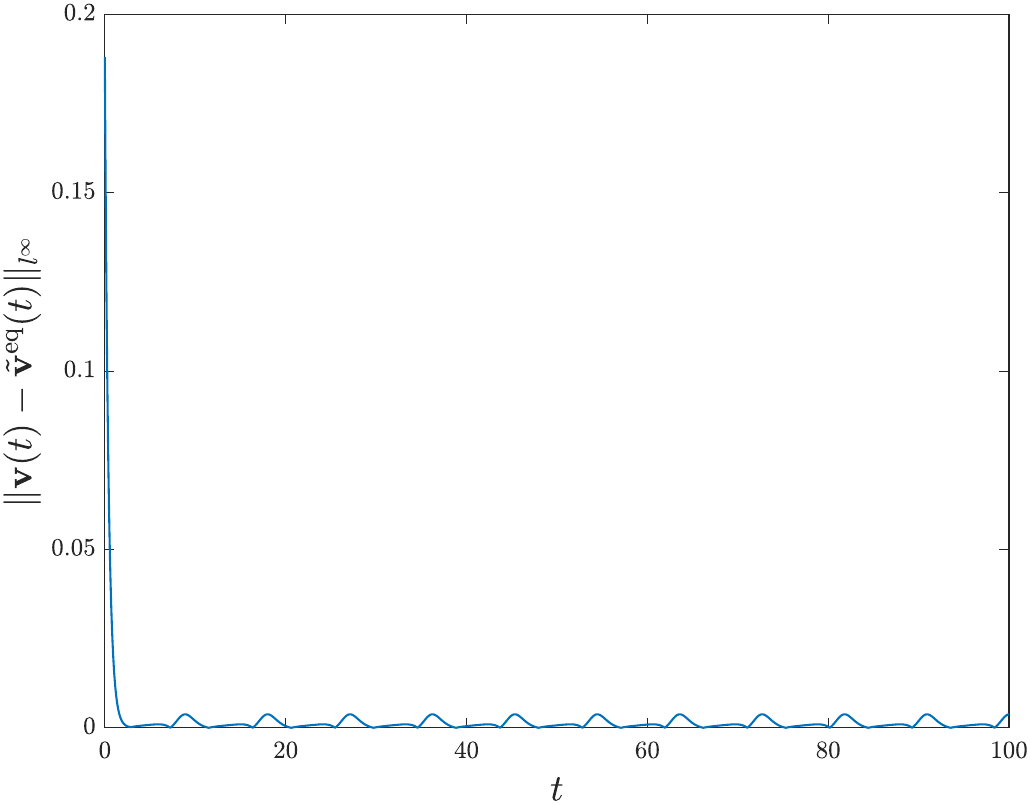}
    \includegraphics[scale=0.3]{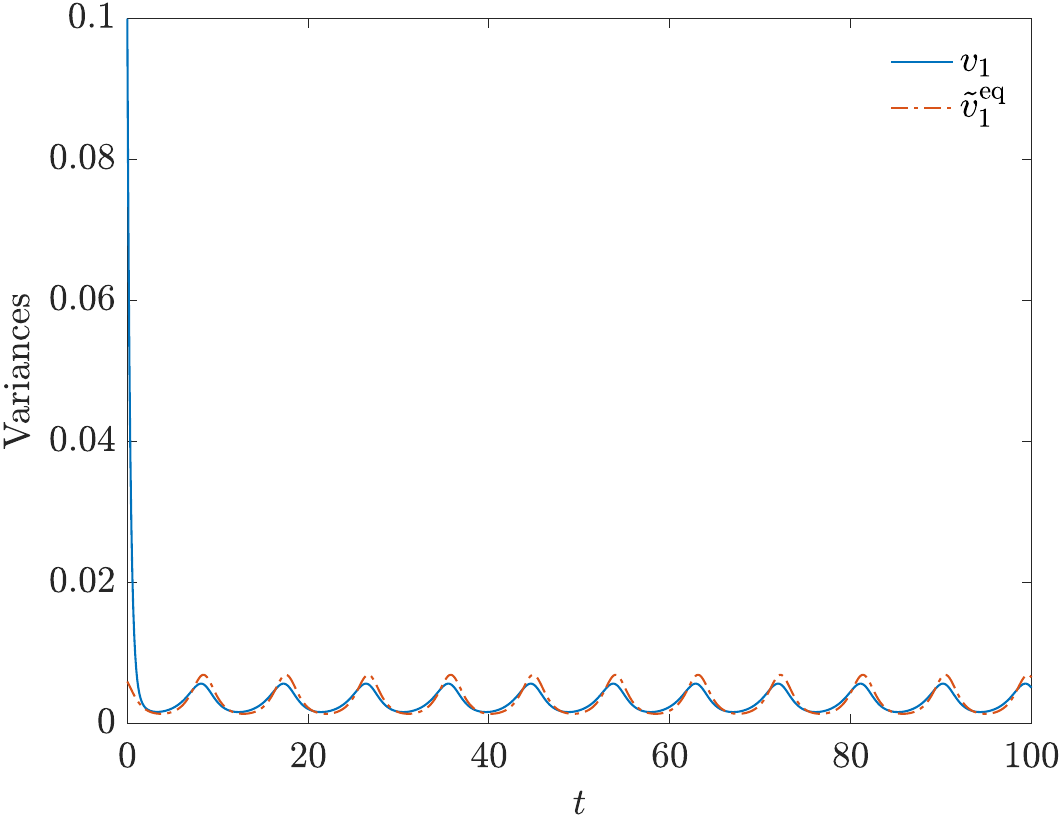}
    \includegraphics[scale=0.3]{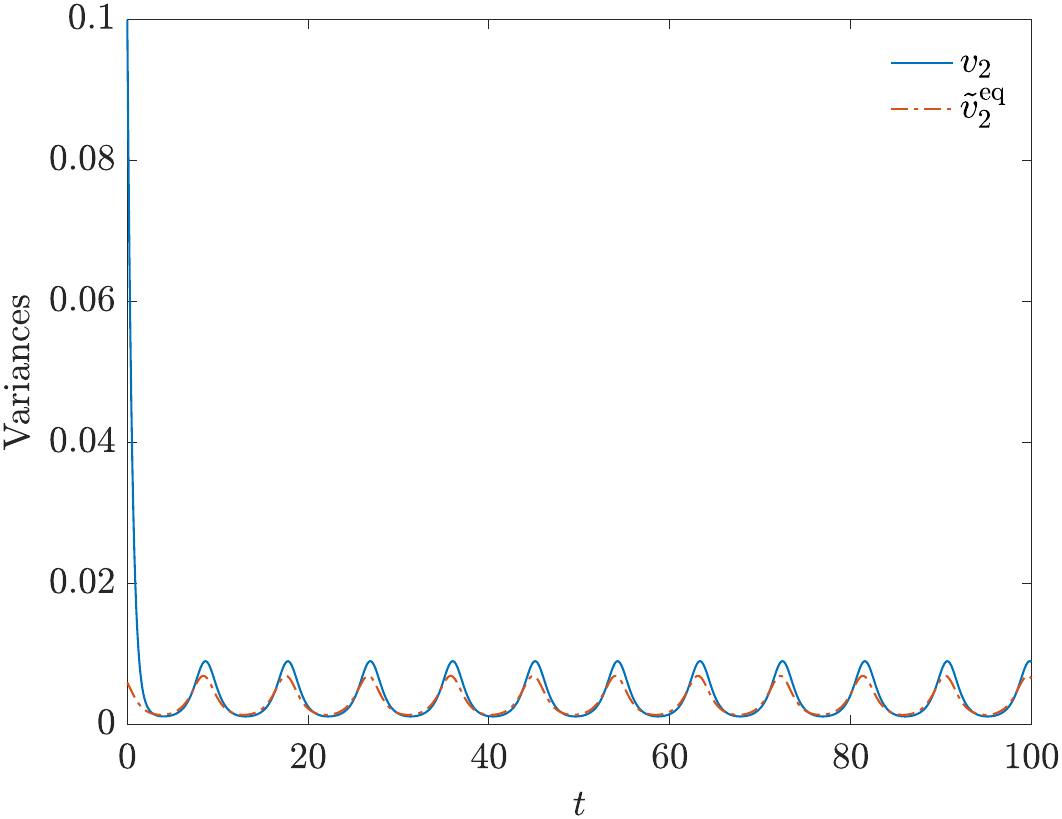}
    \caption{Analysis of the distance between the solution $\mathbf{f}(\x,t)$ to the Fokker--Planck system \eqref{eq:Fokker-Planck} and the local Gamma equilibrium state $\mathbf{f}^\eq(\x,t)$ given by \eqref{eq:equilibrium states}, in terms of the variances $\mathbf{v}(t)$ and $\tilde{\mathbf{v}}^\eq(t)$. Left: temporal evolution of the distance $\|\mathbf{v}(t) - \tilde{\mathbf{v}}^\eq(t)\|_{l^\infty}$. Center: comparison between the dynamics of $v_1(t)$ and $\tilde{v}_1^\eq(t)$. Right: comparison between the dynamics of $v_2(t)$ and $\tilde{v}_2^\eq(t)$. We run the simulations for a time $t \in [0,100]$, with initial conditions $(m_1(0), m_2(0)) = (4,3)$ and $(v_1(0), v_2(0)) = (0.1, 0.1)$.}
\label{fig:distance from rescaled equilibrium}
\end{figure}
%%%%%%%%%%%%%%%%%%%%%%%%%%%%%
%%%%%%%%%%%%%%%%%%%%%%%%%%%%%

In Figure \ref{fig:distance from rescaled equilibrium} we plot the evolution over time of the distance $\|\mathbf{v}(t) - \tilde{\mathbf{v}}^\eq(t)\|_{l^\infty}$ and we compare the dynamics of the two variances, relatively to each population. The simulations show that, while the means $\tilde{\mathbf{m}}^\eq(t)$ of the rescaled populations \eqref{eq:rescaled populations} at equilibrium $\mathbf{f}^\eq(\x,t)$ describe the exact orbits of evolution of the solution $\mathbf{m}(t)$ to the Lotka--Volterra system \eqref{eq:Lotka-Volterra}, their variances $\tilde{\mathbf{v}}^\eq(t)$ are still not able to capture the correct deviations from these orbits. In the next section we propose a possible strategy to stabilize the dynamics of our kinetic model \eqref{eq:Fokker-Planck} and overcome some of the difficulties regarding the investigation of the asymptotic behavior of its solutions.

\begin{remark} \label{remark:equilibrium}
From the previous observations on the means $\mathbf{m}(t)$ and variances $\mathbf{v}(t)$, we may draw similar considerations regarding the relationship between the solution $\mathbf{f}(\x,t)$ and the local equilibrium state $\mathbf{f}^\eq(\x,t)$. In particular, we can infer that the two distributions orbit around a common equilibrium state $\mathbf{f}^*(\x)$ which is the unique uniform-in-time Gamma distribution canceling out the Fokker--Planck fluxes \eqref{eq:equilibrium relations}, given by 
\begin{equation*}
    \begin{split}
        f_1^*(x) & = C_1^* x^{\frac{2\alpha(\chi+1)}{\sigma_1}\frac{m_1^*}{m_2^*} - 1} \exp\left\{-\frac{2}{\sigma_1 m_2^*} (\beta m_2^* + \alpha \chi) x \right\},
        \\[4mm]
        f_2^*(y) & = C_2^* y^{\frac{2 \nu(\theta+1)}{\sigma_2}\frac{m_2^*}{m_1^*} - 1} \exp\left\{-\frac{2}{\sigma_2 m_1^*} (\gamma(\mu - m_1^* + \nu\theta) y \right\},
    \end{split}
\end{equation*}
where $C_1^*$, $C_2^* > 0$ are normalization coefficients depending on the equilibrium means $m_1^*$ and $m_2^*$. One then easily checks that the means and variances of this equilibrium state coincide with the equilibrium points $\mathbf{m}^*$ and $\mathbf{v}^*$, namely 
\begin{align*}
    & \int_{\R_+} x f_1^*(x) \dd x = m_1^*,
    & \int_{\R_+} y f_2^*(y) \dd y = m_2^*, \qquad \\[4mm]
    \int_{\R_+} & (x - m_1^*)^2 f_1^*(x) \dd x = v_1^*,
    & \int_{\R_+} (y - m_2^*)^2 f_2^*(y) \dd y = v_2^*.
\end{align*}
\end{remark}

\begin{remark}
The obtained Fokker--Planck approximation, together with the closure of the moment systems, is inherently limited in capturing the full spectrum of fluctuations near population extinction thresholds, particularly in regimes with finite carrying capacities. Indeed, the modeling assumptions underlying the derivation of the Fokker--Planck system break down in the presence of a finite number of agents, which triggers a regime of rarefied interactions, see \cite{Dobramsyl,Tauber}.
\end{remark}

%%%%%%%%%%%%%%%%%%%%%%%%%%%%%%%%%%%%%%%%%%%%%%%%%%%%%%%%%%%%%%%%%%%
%%%%%%%%%%%%%%%%%%%%%%%%%%%%%%%%%%%%%%%%%%%%%%%%%%%%%%%%%%%%%%%%%%%
%%%%%%%%%%%%%%    SECTION 4: modelingDEL WITH LOGISTIC    %%%%%%%%%%%%%%%
%%%%%%%%%%%%%%%%%%%%%%%%%%%%%%%%%%%%%%%%%%%%%%%%%%%%%%%%%%%%%%%%%%%
%%%%%%%%%%%%%%%%%%%%%%%%%%%%%%%%%%%%%%%%%%%%%%%%%%%%%%%%%%%%%%%%%%%
	
\section{Stabilizing the dynamics with intraspecific interactions}\label{sec:logistic}

The previous section has highlighted the main challenges encountered in studying the relaxation toward equilibrium for the Fokker--Planck system \eqref{eq:Fokker-Planck}, when the underlying macroscopic dynamics (in our case, the Lotka--Volterra ones) exhibit an oscillatory behavior around the equilibrium state. Such a phenomenon appears to render any proof of relaxation toward equilibrium unattainable for the time being. To address this issue, in this section we introduce a correction to our kinetic modeling, switching from a Malthusian-type growth for the preys to a Verhulst-type one, in order to obtain a Lotka--Volterra system with a unique equilibrium state that serves as an attractor for the macroscopic dynamics. This represents another novelty of this work.

Bearing in mind the structure of the Boltzmann equations \eqref{eq:Boltzmann}, we introduce the following system
\begin{equation} \label{eq:Boltzmann logistic}
    \begin{split}
        \frac{\partial f_1(x,t)}{\partial t} &= R_{\chi}^{\alpha}(f_1)(x,t) + Q_{11}(f_1,f_1)(x,t) + Q_{12}(f_1,f_2)(x,t),
        \\[2mm]
        \frac{\partial f_2(y,t)}{\partial t} &= R_{\theta}^{\nu}(f_2)(y,t) + Q_{21}(f_2,f_1)(y,t),
    \end{split}
\end{equation}
which accounts for both interspecific predator--preys interactions as well as intraspecific competition among the preys. More precisely, the bilinear operators $Q_{12}(f_1,f_2)$ and $Q_{21}(f_2,f_1)$ are now used to describe any binary interaction between individuals of the same or of different species, and the operator $Q_{11}(f_1,f_1)$ models in particular the intraspecific interactions occurring inside the population of preys. This means that at the microscopic level, the preys are involved in the interactions \eqref{eq:microscopic Boltzmann} with the predators, and in new intraspecific binary interactions given by
\begin{equation} \label{eq:microscopic Boltzmann logistic}
        x''' = x - \tilde{\Phi}(x_*) x + \sqrt{x} \mathds{1}(x \geq s_0/2) \eta_1(x_*),
\end{equation}
where the variable $x_* \in \R_+$ refers to the size of the interacting group within the species of preys, described by $f_1(x_*,t)$, and compared to \eqref{eq:microscopic Boltzmann} we now solely consider the more interesting case $p = 1/2$ which leads to thin-tails equilibrium states. Notice that we have once again fixed $\bar{x} = 1$ for simplicity and we still need to reduce the range of values taken by $x$ by using the indicator function $\mathds{1}(x \geq s_0/2)$, which ensures the positivity of the post-interactions $x'''$. For consistency with our previous choices \eqref{eq:deterministic effects}--\eqref{eq:random effects} on the form of the deterministic and random effects involved in the interactions, we choose once more
\begin{equation*}
    \tilde{\Phi}(x_*) = \frac{\alpha}{K} \frac{x_*}{1+x_*},
\end{equation*}
where the new parameter $K > 0$ denotes the carrying capacity of the preys, and we consider $\eta_1(x_*)$ to be an independent random variable of zero mean and bounded variance, given explicitly by
\begin{equation*}
    \langle\eta_1^2(x_*)\rangle = \sigma_1 \frac{x_*}{1+x_*}.
\end{equation*}
In particular, we assume that the ratio $\alpha/K$ satisfies the constraint $\displaystyle \frac{\alpha}{K} < 1$.

Then, similarly to the previous section, we consider the weak formulation of system \eqref{eq:Boltzmann logistic} which, for all smooth functions $\varphi(x)$ and $\psi(y)$, reads as
\begin{equation} \label{eq:Boltzmann logistic weak}
    \begin{split}
        \frac{\dd}{\dd t}\int_{\R_+ }\varphi(x)f_1(x,t) \dd x &= \int_{\R_+\times \R_+} \kappa_1(x,x_*) \big\langle\varphi(x''')-\varphi(x)\big\rangle f_1(x,t) f_1(x_*,t) \dd x \dd x_* \\[2mm]
        & \hspace*{1cm} +\int_{\R_+ \times \R_+} \big(\varphi(x'')-\varphi(x)\big) f_1(x,t) g(z,t) \dd z \dd x \\[2mm]
        & \hspace*{2cm} +\int_{\R_+\times \R_+} \kappa_1(x,y) \big\langle\varphi(x')-\varphi(x)\big\rangle f_1(x,t)f_2(y,t) \dd x \dd y, \\[4mm]
        \frac{\dd}{\dd t} \int_{\R_+} \psi(y) f_2(y,t) \dd y 
        &= \int_{\R_+\times \R_+} \big(\psi(y'')-\psi(y)\big) f_2(y,t) h(z,t) \dd z \dd y \\[2mm]
        & \hspace*{2cm} +\int_{\R_+\times\R_+} \kappa_2(x,y) \big\langle\psi(y')-\psi(y)\big\rangle f_1(x,t)f_2(y,t) \dd x \dd y,
    \end{split}
\end{equation}
where we have used the same function $\kappa_1(x,x_*) = \kappa(x_*) = 1+x_*$ defined in \eqref{eq:interaction frequencies} to model the frequency of contacts \eqref{eq:microscopic Boltzmann logistic} among preys, and the integral terms on the right-hand sides of the equations account for the balance of all intraspecific \eqref{eq:microscopic Boltzmann logistic} and interspecific \eqref{eq:microscopic Boltzmann} microscopic interactions, as well as of the exchanges with the environment \eqref{eq:microscopic redistribution}, recalling the notations $g(z,t)$ and $h(z,t)$ for the distributions of its resources, that have prescribed means given by \eqref{eq:mean of environment}.

In this new setting, choosing $\varphi(x) = x$ and $\psi(y) = y$ in \eqref{eq:Boltzmann logistic weak}, we find that the temporal evolution of the mean values $m_1(t)$ and $m_2(t)$ is now described by the system
\begin{equation} \label{eq:Lotka-Volterra logistic}
    \begin{split}
        \frac{\dd m_1(t)}{\dd t} &= \alpha \left( 1 - \frac{m_1(t)}{K}\right) m_1(t) - \beta m_1(t)m_2(t),
        \\[2mm]
        \frac{\dd m_2(t)}{\dd t} &= -\delta m_2(t) + \gamma m_1(t)m_2(t),
    \end{split}
\end{equation}
which are the well-known competitive predator--prey Lotka--Volterra equations, involving a logistic growth for the population of preys. It is worth pointing out that this new system differs from the one \eqref{eq:Lotka-Volterra} derived in the previous section only in the first equation, as we solely consider intraspecific interactions among preys. Furthermore, we highlight that system \eqref{eq:Lotka-Volterra logistic} is characterized by a unique stable equilibrium point of coexistence given by
\begin{equation} \label{eq:m^infty}
    \mathbf{m}^\infty = \left( \frac{\delta}{\gamma}, \frac{\alpha (\gamma K - \delta)}{\beta \gamma K} \right),
\end{equation}
which is biologically meaningful provided that $\gamma K - \delta > 0$. In particular, under this condition it is possible to prove that any  solution $\mathbf{m}(t)$ to \eqref{eq:Lotka-Volterra logistic} converges asymptotically in time to the equilibrium $\mathbf{m}^\infty$, see e.g. \cite{ChaJabRao}.

We proceed to derive the corresponding reduced model of Fokker--Planck-type, by considering the following quasi-invariant regime of the parameters
\begin{equation*}
    \alpha \to \eps \alpha, \qquad \beta \to \eps \beta, \qquad \nu \to \eps \nu, \qquad \gamma \to \eps \gamma, \qquad 
    \sigma_1 \to \eps \sigma_1, \qquad \sigma_2 \to \eps \sigma_2,
\end{equation*}
for some $\eps \ll 1$, where the temporal variable is rescaled accordingly as $\tau = \eps t$. Since the computations do not deviate substantially from the previous ones, we skip the details of such derivation. When $\eps \to 0$, we find that the evolution of the limit distribution functions $f_1(x,t)$ and $f_2(y,t)$ is described by the Fokker--Planck system
\begin{equation} \label{eq:Fokker-Planck logistic}
    \begin{split}
        \frac{\partial f_1(x,t)}{\partial t} &= \frac{\sigma_1 (m_1(t) + m_2(t))}{2} \frac{\partial^2}{\partial x^2} \left[ x f_1(x,t) \right] \\[2mm]
        & \hspace*{2.3cm} + \frac{\partial}{\partial x} \left[ \left( \left( \beta m_2(t) + \frac{\alpha}{K} m_1(t) + \alpha \chi \right) x -\alpha(\chi+1) m_1(t) \right) f_1(x,t) \right], \\[6mm]
        \frac{\partial f_2(y,t)}{\partial t} &= \frac{\sigma_2 m_1(t)}{2} \frac{\partial^2}{\partial y^2} \left[y f_2(y,t)\right] + \frac{\partial}{\partial y} \big( (\gamma(\mu - m_1(t)) y + \nu\theta y -\nu(\theta+1) m_2(t)) f_2(y,t) \big),
    \end{split}
\end{equation}
complemented by the no-flux boundary conditions
\begin{equation*}
    \begin{split}
        & \frac{\sigma_1 (m_1(t) + m_2(t))}{2} \frac{\partial}{\partial x} \left[x f_1(x,t)\right] \\[2mm]
        & \hspace*{1.7cm} + \left( \left( \beta m_2(t) + \frac{\alpha}{K} m_1(t) + \alpha \chi \right) x -\alpha(\chi+1) m_1(t) \right) f_1(x,t) \Big|_{x=0, +\infty} = 0, \\[2mm]
        & x f_1(x,t) \Big|_{x=0, +\infty} = 0, \\[8mm]
        & \frac{\sigma_2 m_1(t)}{2} \frac{\partial}{\partial y} \left[y f_2(y,t)\right] + (\gamma(\mu - m_1(t)) y + \nu\theta y - \nu(\theta+1) m_2(t)) f_2(y,t) \Big|_{y=0, +\infty} = 0, \\[2mm]
        & y f_2(y,t) \Big|_{y=0, +\infty} = 0. 
\end{split}
\end{equation*}
Starting from \eqref{eq:Fokker-Planck logistic}, we are now able to determine closed evolution equations for the moments of the system, in particular for the associated means and variances of its solutions $\mathbf{f}(\x,t)$. For the means $m_1(t)$ and $m_2(t)$ we find a governing system that is obviously consistent with the one determined from the Boltzmann dynamics, given by the competitive Lotka--Volterra equations \eqref{eq:Lotka-Volterra logistic}. The evolution of the variances $v_1(t)$ and $v_2(t)$ is then prescribed by 
\begin{equation} \label{eq:variance FP logistic}
    \begin{split}
        \frac{\dd v_1(t)}{\dd t} &= -2 \left( \beta m_2(t) + \frac{\alpha}{K} m_1(t) + \alpha \chi \right) v_1(t) + \sigma_1 m_1(t)(m_1(t) + m_2(t)), \\[4mm]
        \frac{\dd v_2(t)}{\dd t} &= -2 \big( \gamma(\mu - m_1(t)) + \nu\theta \big) v_2(t) + \sigma_2 m_1(t) m_2(t),
    \end{split}
\end{equation}
and one can infer global-in-time a priori estimates similar to the ones \eqref{eq:bounds on variances} obtained in the previous section.

Moreover, we can carry out the same computations as before to determine the new local equilibrium distribution states $\mathbf{f}^\eq(\x,t) = (f_1^\eq(x,t),f_2^\eq(y,t))$ which now read
\begin{equation} \label{eq:equilibrium states logistic}
    \begin{split}
        f_1^\eq(x,t) &= C_1(t) x^{\frac{2\alpha(\chi+1)}{\sigma_1}\frac{m_1(t)}{m_1(t) + m_2(t)} - 1} \exp\left\{-\frac{2}{\sigma_1 (m_1(t) + m_2(t))} \left(\beta m_2(t) + \frac{\alpha}{K} m_1(t) + \alpha \chi \right) x \right\}, \\[4mm]
        f_2^\eq(y,t) &= C_2(t) y^{\frac{2 \nu(\theta+1)}{\sigma_2}\frac{m_2(t)}{m_1(t)} - 1} \exp\left\{-\frac{2}{\sigma_2 m_1(t)} (\gamma(\mu - m_1(t)) + \nu\theta) y \right\},
    \end{split}
\end{equation}
for some time-dependent normalization constants $C_1(t), C_2(t) > 0$. The above functions are Gamma-type distributions if we impose (once the parameters $\alpha$, $K$, $\beta$, $\delta$, $\gamma$, $\nu > 0$ and $\mu \geq 1$ of the competitive Lotka--Volterra system \eqref{eq:Lotka-Volterra logistic} have been fixed) that $\chi,\theta \in \R$ satisfy the conditions
\begin{equation} \label{eq:admissible parameters logistic}
    \chi > \max\left\{ -1,\ -\frac{\beta}{\alpha}\underline{c}_2 - \frac{1}{K}\underline{c}_1 \right\}, \qquad \theta > \max\left\{ -1,\ -\frac{\gamma}{\nu} (\mu - \overline{c}_1) \right\},
    %\begin{array}{llll}
    %    \displaystyle \chi > - \max\left\{ 1,\ \frac{\beta}{\alpha}\underline{c}_2 + \frac{1}{K} \underline{c}_1 \right\} & \textrm{and} & \displaystyle \theta > 0 & (\textrm{if} \ \mu \geq \overline{c}_1), \\[6mm]
    %    \displaystyle \chi > - \max\left\{ 1,\ \frac{\beta}{\alpha}\underline{c}_2 + \frac{1}{K} \underline{c}_1 \right\} & \textrm{and} & \displaystyle \theta > \frac{\gamma}{\nu} (\overline{c}_1 - \mu) & (\textrm{if} \ \underline{c}_1 \leq \mu < \overline{c}_1).
    %\end{array}
\end{equation}
for some constants $\underline{c}_1$, $\overline{c}_1$, $\underline{c}_2 > 0$ bounding from below and from above the means $m_1(t)$ and $m_2(t)$, defined similarly to \eqref{eq:bounds on means}. Notice in particular that the introduction of the intraspecific competition among preys allows to lessen the restrictions on the range of admissible values for $\chi$, meaning that the local equilibrium $\mathbf{f}^\eq(\x,t)$ is well-defined as a Gamma vector distribution and the corresponding variances $v_1(t)$ and $v_2(t)$ remain bounded, for values of $\chi$ ranging in an interval wider than the one determined in the previous setting.

By rewriting the equilibrium distributions \eqref{eq:equilibrium states logistic} as
\begin{equation*} %\label{eq:qe_logistic}
    f_1^\eq(x,t) = C_1(t) x^{a_1(t) - 1} e^{- b_1(t) x}, \qquad f_2^\eq(y,t) = C_2(t) y^{a_2(t) - 1} e^{- b_2(t) y},
\end{equation*}
where all the exponents
\begin{equation*} %\label{parameters_logistic}
    \begin{split}
        & a_1(t) = \frac{2\alpha(\chi+1)m_1(t)}{\sigma_1 (m_1(t)+m_2(t))}, \qquad b_1(t) = \frac{2(\beta m_2(t) + \frac{\alpha}{K} m_1(t) + \alpha \chi)}{\sigma_1(m_1(t)+m_2(t))}, \\[4mm]
        & a_2(t) = \frac{2\nu(\theta+1)m_2(t)}{\sigma_2 m_1(t)}, \qquad \quad\ b_2(t) = \frac{2(\gamma(\mu-m_1(t))+\nu\theta)}{\sigma_2 m_1(t)},
    \end{split}
\end{equation*}
are positive thanks to the above assumptions on the parameters $\chi$ and $\theta$, we deduce explicit expressions for the equilibrium means $\mathbf{m}^\eq(t)$, which are given by
\begin{equation} \label{eq:Gamma means logistic}
    \begin{split}
        m_1^\eq(t) &= \frac{\alpha(\chi+1)m_1(t)}{\beta m_2(t) + \frac{\alpha}{K} m_1(t) + \alpha \chi}, \\[4mm]
        m_2^\eq(t) &= \frac{\nu(\theta + 1) m_2(t)}{\gamma(\mu - m_1(t)) + \nu\theta},
    \end{split}
\end{equation}
and for the equilibrium variances $\mathbf{v}^\eq(t)$, which read
\begin{equation} \label{eq:Gamma variances logistic}
    \begin{split}
        v_1^\eq(t) &= \frac{\sigma_1 \alpha(\chi+1) m_1(t) (m_1(t)+m_2(t))}{2\left(\beta m_2(t) + \frac{\alpha}{K} m_1(t) + \alpha \chi \right)^2} = \frac{\sigma_1 m_1^\eq(t) (m_1(t)+m_2(t))}{2\left(\beta m_2(t) + \frac{\alpha}{K} m_1(t) + \alpha \chi \right)}, \\[4mm]
        v_2^\eq(t) &= \frac{\sigma_2 \nu(\theta+1) m_1(t) m_2(t)}{2(\gamma(\mu-m_1(t)) + \nu \theta)^2} = \frac{\sigma_2 m_1(t) m_2^\eq(t)}{2(\gamma(\mu-m_1(t)) + \nu \theta)}. 
    \end{split}
\end{equation}
We can thus investigate the asymptotic behavior of the first two moments of the distributions of preys and predators, in order to evaluate the distance of the solutions $\mathbf{f}(\x,t)$ to system \eqref{eq:Fokker-Planck logistic} from the associated local equilibrium states $\mathbf{f}^\eq(\x,t)$ given by \eqref{eq:equilibrium states logistic}. It is easy to check that the vector $\mathbf{m}^\infty$ defined by \eqref{eq:m^infty} is the unique equilibrium point of both $\mathbf{m}(t)$ and $\mathbf{m}^\eq(t)$. Similarly, one can observe that the unique equilibrium point $\mathbf{v}^\infty$, common to $\mathbf{v}(t)$ and $\mathbf{v}^\eq(t)$, is the one obtained by evaluating at $\mathbf{m}^\infty$ the couples $(v_1(t),v_2(t))$
%\begin{equation*}
%    \begin{split}
%        v_1(t) = \frac{\sigma_1 m_1(t) (m_1(t)+m_2(t)}{2 \left(\beta m_2(t) + \frac{\alpha}{K} + \alpha \chi \right)}, \qquad v_2(t) = \frac{\sigma_2 m_1(t) m_2(t)}{2 (\gamma (\mu - m_1(t)) + \nu\theta)},
%    \end{split}
%\end{equation*}
canceling out the right-hand sides of equations \eqref{eq:variance FP logistic}. In particular, we recover
\begin{equation} \label{eq:v^infty}
    \mathbf{v}^\infty = \left(\frac{\sigma_1 m_1^\infty (m_1^\infty + m_2^\infty)}{2(\beta m_2^\infty  + \frac{\alpha}{K}m_1^\infty + \alpha \chi)}, \frac{\sigma_2 m_1^\infty m_2^\infty}{2 \left( \gamma (\mu -m_1^\infty) + \nu \theta \right)}\right),
\end{equation}
which coincides with the limit for $t \to +\infty$ of the right-hand sides (specifically, the second equalities) of relations \eqref{eq:Gamma variances logistic}, as soon as one proves that both $\mathbf{m}(t) \underset{t \to +\infty}{\longrightarrow} \mathbf{m}^\infty$ and $\mathbf{m}^\eq(t) \underset{t \to +\infty}{\longrightarrow} \mathbf{m}^\infty$. The first convergence is obvious, since $\mathbf{m}^\infty$ is the unique globally asymptotically stable equilibrium of the competitive Lotka--Volterra system \eqref{eq:Lotka-Volterra logistic}, while the second can be recovered by taking the limit $t \to +\infty$ of relations \eqref{eq:Gamma means logistic} and performing simple algebraic manipulations.

%%%%%%%%%%%%%%%%%%%%%%%%%%%%%
%%%%%%%%%%%%%%%%%%%%%%%%%%%%%
\begin{figure}
\centering
    \includegraphics[scale = 0.4]{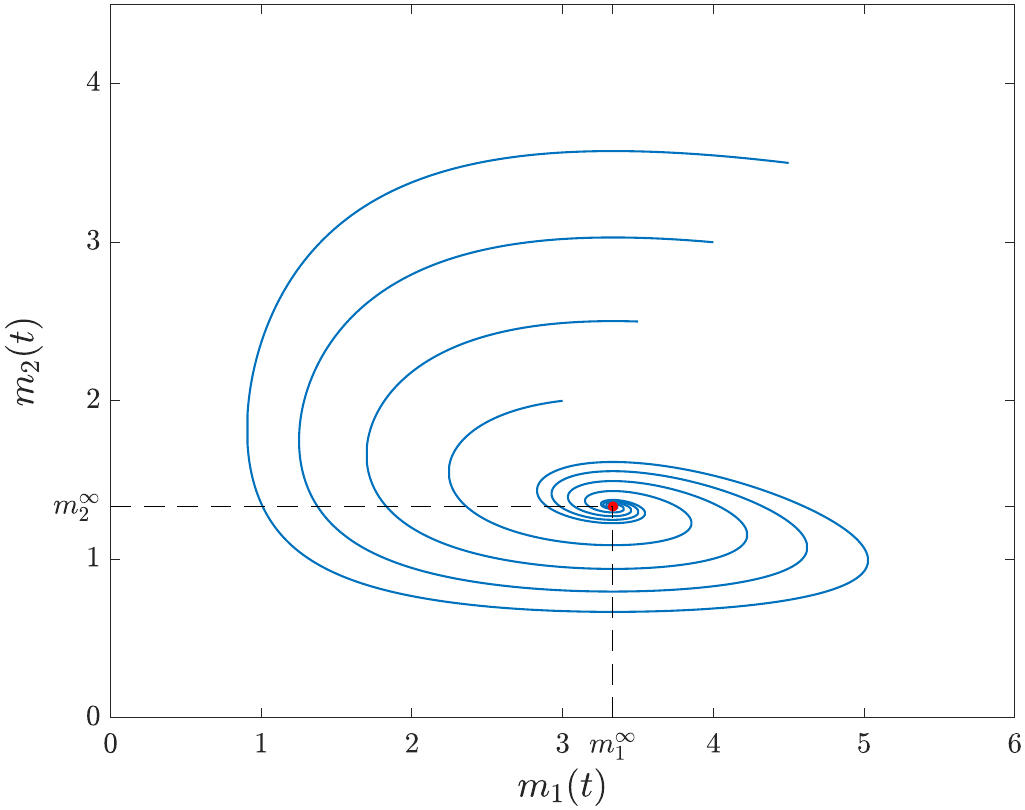}
    \includegraphics[scale = 0.4]{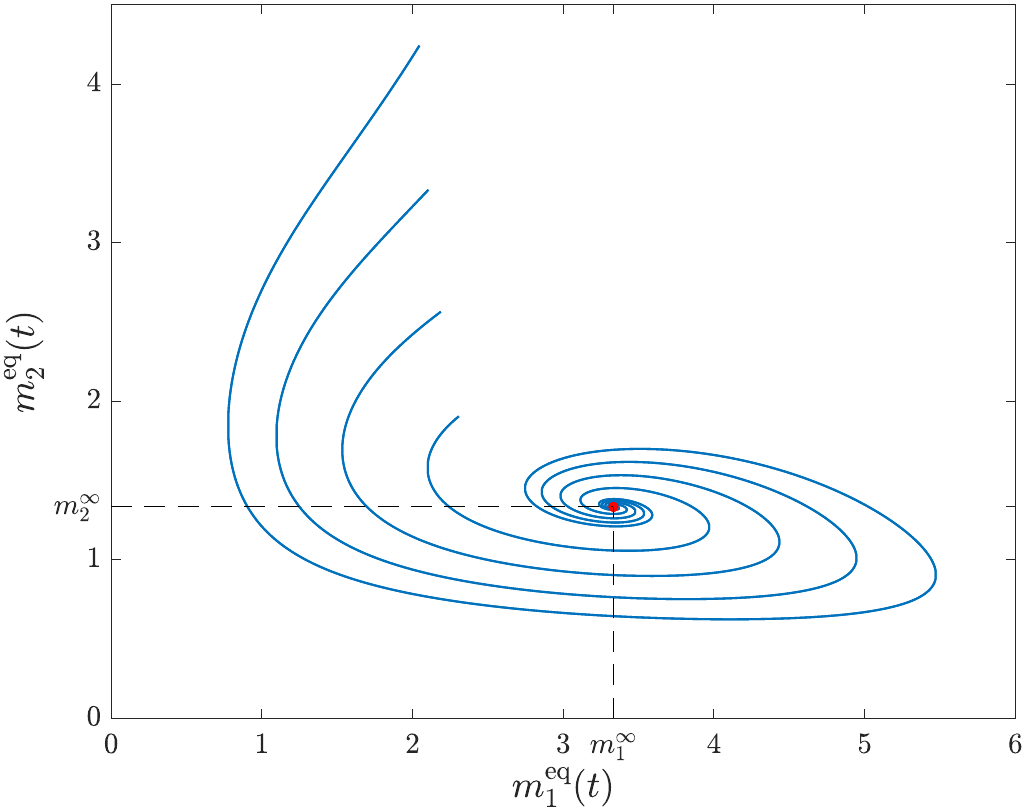} \\
    \includegraphics[scale = 0.4]{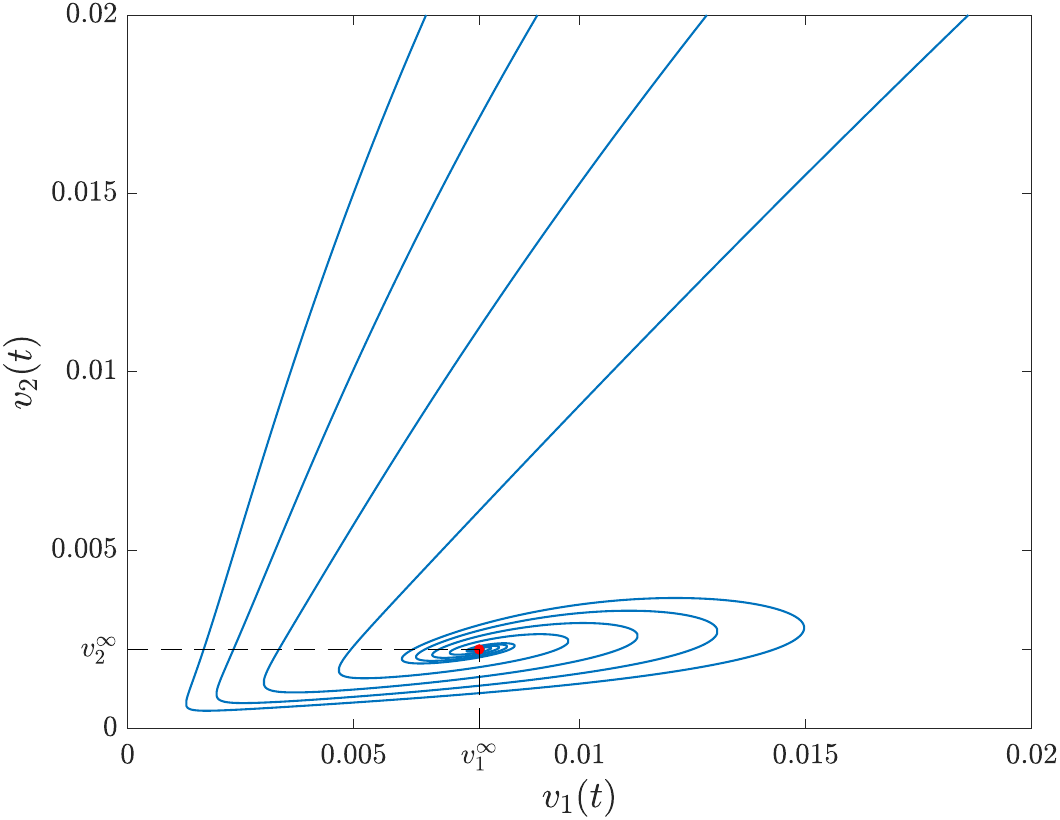} 
    \includegraphics[scale = 0.4]{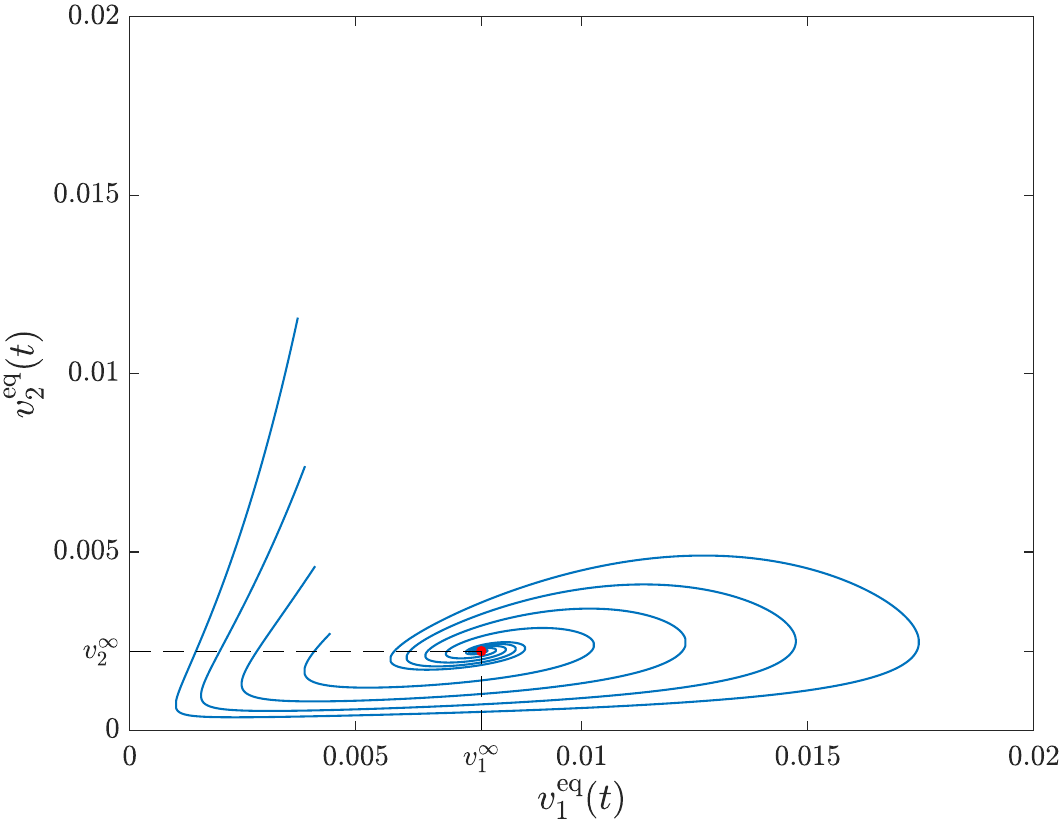}
    \caption{Top row: orbits of evolution of the solutions $\mathbf{m}(t)$ to the logistic Lotka--Volterra model \eqref{eq:Lotka-Volterra logistic} and of the means $\mathbf{m}^\eq(t)$ of the local equilibrium distributions defined in \eqref{eq:equilibrium states logistic}. In red we plot the equilibrium point $\mathbf{m}^\infty$. Bottom row: orbits of evolution of the variances $\mathbf{v}(t)$ solving system \eqref{eq:variance FP logistic} and of the variances $\mathbf{v}^\eq(t)$ computed from the corresponding equilibrium state \eqref{eq:equilibrium states logistic}. In red we plot the point $\mathbf{v}^\infty$. As initial data we considered $(m_1(0),m_2(0)) = (4,3)$ and $v_1(0) = v_2(0) = 0.1$.}
    \label{fig:orbits logistic}
\end{figure}
%%%%%%%%%%%%%%%%%%%%%%%%%%%%%
%%%%%%%%%%%%%%%%%%%%%%%%%%%%%

%This shows the main difference with respect to the framework studied in the previous section. Indeed, if we investigate the asymptotic behavior of the first two moments of the distributions of preys and predators, in order to evaluate the distance of the solutions $\mathbf{f}(\x,t)$ to system \eqref{eq:Fokker-Planck logistic} from the associated local equilibrium states $\mathbf{f}^\eq(\x,t)$ given by \eqref{eq:equilibrium states logistic}, 
In particular, we can prove that the means $\mathbf{m}(t)$ and $\mathbf{m}^\eq(t)$ satisfy the relative convergence estimate
\begin{equation} \label{eq:main estimate logistic}
    \begin{split}
        \| \mathbf{m}(t) - \mathbf{m}^\eq(t) \|_{l^\infty} &= \frac{\beta \left| m_2(t) - \frac{\alpha}{\beta} \left(1 - \frac{m_1(t)}{K} \right) \right|}{\beta m_2(t) + \frac{\alpha}{K} m_1(t) + \alpha \chi} m_1(t) + \frac{\gamma \left| m_1(t) - \frac{\delta}{\gamma} \right|}{\gamma (\mu - m_1(t)) + \nu \theta} m_2(t) \\[6mm]
        & \leq \frac{\beta K + \alpha}{\beta K} \max \left\{ \frac{\beta \overline{c}_1}{\beta \underline{c}_2 + \frac{\alpha}{K} \underline{c}_1 + \alpha \chi}, \frac{\gamma \overline{c}_2}{\gamma(\mu - \overline{c}_1) + \nu\theta} \right\} \| \mathbf{m}(t) - \mathbf{m}^{\infty} \|_{l^\infty} \underset{t \to +\infty}{\longrightarrow} 0.
    \end{split}
\end{equation}
Additionally, the equilibrium point $\mathbf{v}^\infty$ is also globally asymptotically stable for system \eqref{eq:variance FP logistic}, because the latter can be recast as
\begin{equation*}
    \begin{split}
        \frac{\dd}{\dd t} (v_1(t) - v_1^\infty) &= -2 \left( \beta m_2(t) + \frac{\alpha}{K} m_1(t) + \alpha \chi \right) (v_1(t) - v_1^\infty) \\[2mm]
        & \hspace*{2cm} + 2 \left( \beta m_2(t) + \frac{\alpha}{K} m_1(t) + \alpha \chi \right) \left( \frac{\sigma_1 m_1(t)(m_1(t) + m_2(t))}{2 \left( \beta m_2(t) + \frac{\alpha}{K} m_1(t) + \alpha \chi \right)} - v_1^\infty \right), \\[4mm]
        \frac{\dd}{\dd t} (v_2(t) - v_2^\infty) &= -2 \big( \gamma(\mu - m_1(t)) + \nu\theta \big) (v_2(t) - v_2^\infty) \\[2mm]
        & \hspace*{2cm} + 2 \big( \gamma(\mu - m_1(t)) + \nu\theta \big) \left( \frac{\sigma_2 m_1(t) m_2(t)}{2 \big( \gamma(\mu - m_1(t)) + \nu\theta \big)} - v_2^\infty \right),
    \end{split}
\end{equation*}
from which one deduces, via Gronwall's inequality, the bounds
\begin{equation*}
    \begin{split}
        |v_1(t) - v_1^\infty| &\leq |v_1(0) - v_1^\infty| e^{-2 \zeta_1 t} \\[2mm] 
        & \hspace*{2cm} + \frac{\beta \overline{c}_2 + \frac{\alpha}{K} \overline{c}_1 + \alpha |\chi|}{\zeta_1} \left| \frac{\sigma_1 m_1(t)(m_1(t) + m_2(t))}{2 \left( \beta m_2(t) + \frac{\alpha}{K} m_1(t) + \alpha \chi \right)} - v_1^\infty \right| \underset{t \to +\infty}{\longrightarrow} 0, \\[4mm]
        |v_2(t) - v_2^\infty| &\leq |v_2(0) - v_2^\infty| e^{-2 \zeta_2 t} \\
        & \hspace*{3.2cm} + \frac{\gamma(\mu + \overline{c}_1) + \nu |\theta|}{\zeta_2} \left| \frac{\sigma_2 m_1(t) m_2(t)}{2 \big( \gamma(\mu - m_1(t)) + \nu\theta \big)} - v_2^\infty \right| \underset{t \to +\infty}{\longrightarrow} 0,
    \end{split}
\end{equation*}
where $\zeta_1, \zeta_2 > 0$ denotes two constants such that $\beta \underline{c}_2 + \frac{\alpha}{K} \underline{c}_1 + \alpha \chi \geq \zeta_1$ and $\gamma(\mu - \overline{c}_1) + \nu\theta \geq \zeta_2$, whose existence is guaranteed by conditions \eqref{eq:admissible parameters logistic} on $\chi$ and $\theta$. The above convergences come from the definition of $\mathbf{v}^\infty$, by using that $\mathbf{m}(t) \underset{t \to +\infty}{\longrightarrow} \mathbf{m}^\infty$, and they imply its global asymptotic stability $\| \mathbf{v}(t) - \mathbf{v}^\infty \|_{l^\infty} \underset{t \to +\infty}{\longrightarrow} 0$, as desired. Owing to the fact that $\mathbf{v}^\eq(t) \underset{t \to +\infty}{\longrightarrow} \mathbf{v}^\infty$, which can be proved by taking the limit $t \to +\infty$ into relations \eqref{eq:Gamma variances logistic}, we conclude that also the relaxation $\| \mathbf{v}(t) - \mathbf{v}^\eq(t) \|_{l^\infty} \underset{t \to +\infty}{\longrightarrow} 0$ hols true. In Figure \ref{fig:orbits logistic} we display the different orbits of evolution for the means and variances of the solution $\mathbf{f}(\x,t)$ and of the corresponding local Gamma equilibrium state $\mathbf{f}^\eq(\x,t)$.

These considerations show the main difference with respect to the framework studied in the previous section. Indeed, the introduction of intraspecific interactions leads as expected to more stable dynamics for the evolution of the means, allowing to recover a relaxation in time of $\mathbf{m}(t)$ toward $\mathbf{m}^\eq(t)$. This is due to the fact that the logistic term makes the macroscopic equilibrium point $\mathbf{m}^\infty$ an attractor for the system, thus identifying a unique orbit of evolution for the means $\mathbf{m}(t)$. In particular, $\mathbf{m}(t)$ and $\mathbf{m}^\eq(t)$ are still localized on different orbits, but they now converge to the common limit $\mathbf{m}^\infty$, meaning that they must be close in their large-time asymptotics, as expressed by estimate \eqref{eq:main estimate logistic}. This stabilization of the dynamics for the moments of the solutions $\mathbf{f}(\x,t)$ to \eqref{eq:Lotka-Volterra logistic} is further confirmed by studying the asymptotic behavior of their variances $\mathbf{v}(t)$ which display a relaxation toward the equilibrium point $\mathbf{v}^\infty$, in common with the variances $\mathbf{v}^\eq(t)$ of the associated local equilibrium $\mathbf{f}^\eq(\x,t)$ defined by \eqref{eq:equilibrium states logistic}, resulting in the relative convergence of $\mathbf{v}(t)$ toward $\mathbf{v}^\eq(t)$. This means that $\mathbf{m}^\eq(t)$ and $\mathbf{v}^\eq(t)$ provide good predictions for the means $\mathbf{m}(t)$ and variances $\mathbf{v}(t)$ respectively, in the longtime asymptotics of the Fokker--Planck system \eqref{eq:Fokker-Planck logistic}.

\begin{figure}
\centering  
    \includegraphics[scale=0.4]{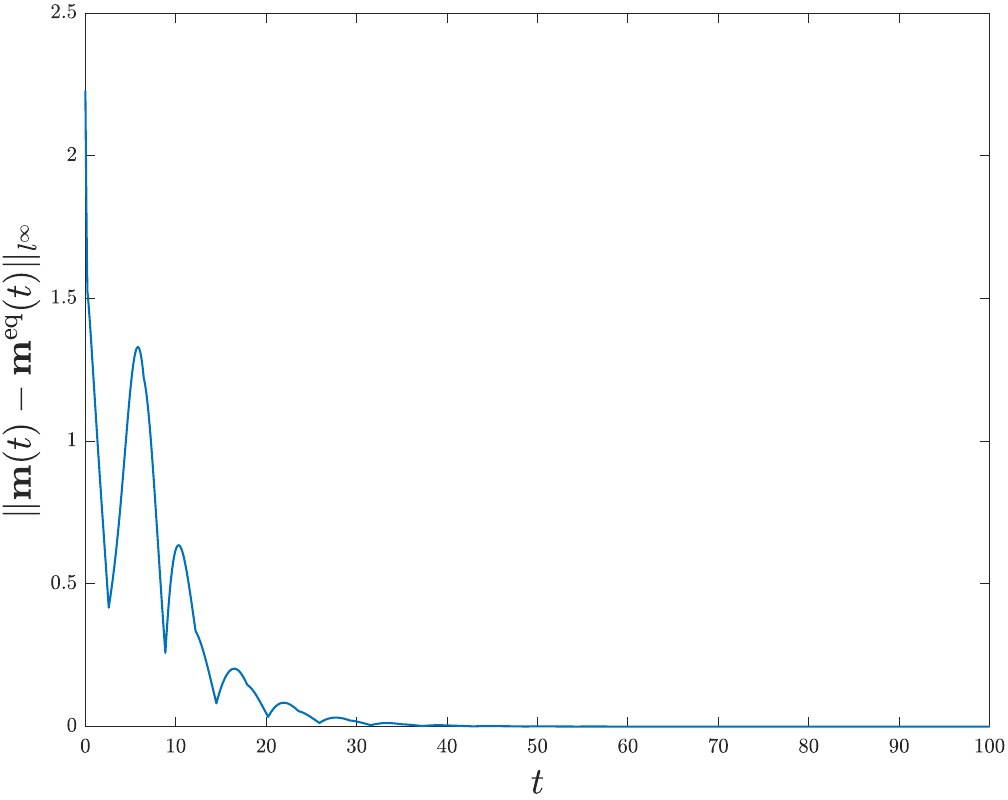} \hspace*{5mm}
    \includegraphics[scale=0.4]{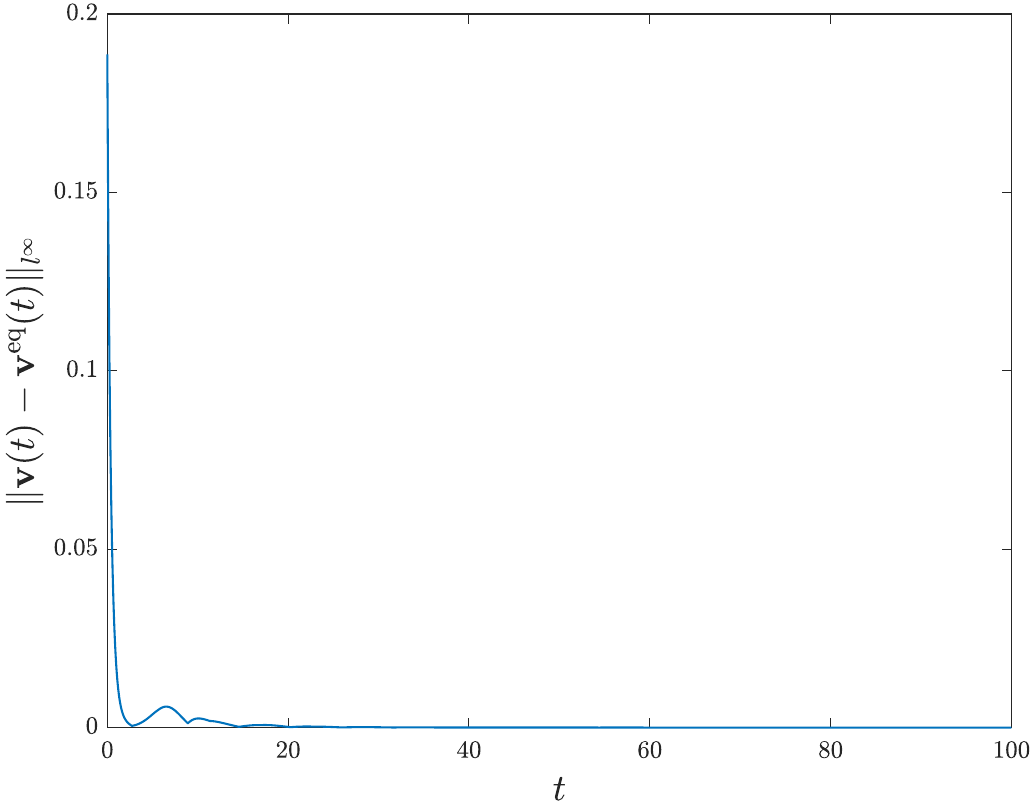}
    \caption{Analysis of the distance between the the solution $\mathbf{f}(\x,t)$ to the competitive Fokker--Planck system \eqref{eq:Fokker-Planck logistic} and local Gamma equilibrium state $\mathbf{f}^\eq(\x,t)$ given by \eqref{eq:equilibrium states logistic}, in terms of their moments of order one and two. Evolution over time of the distances between the respective means (left) and variances (right). We run the simulation for a time $t \in [0,100]$, with initial conditions $(m_1(0), m_2(0)) = (4,3)$ and $v_1(0) = v_2(0) = 0.1$.}
\label{fig:distance from equilibrium logistic}
\end{figure}

In Figure \ref{fig:distance from equilibrium logistic} we show the results of our numerical simulations, evaluating the distance between the solutions $\mathbf{v}(t)$ to system \eqref{eq:variance FP logistic} and the variances $\mathbf{v}^\eq(t)$ computed from the local equilibrium states \eqref{eq:equilibrium states logistic}. To recover the evolution of $\mathbf{v}(t)$, we simulate in particular the coupling of systems \eqref{eq:Lotka-Volterra logistic} and \eqref{eq:variance FP logistic}, for the means and the variances respectively. The parameters are given in Table \ref{tab:parameters}, and we also fix $K = 10$. On the left, we can see that the means $\mathbf{m}(t)$ converge to the corresponding local equilibrium $\mathbf{m}^\eq(t)$, in agreement with the analytical result \eqref{eq:main estimate logistic}. On the right, the numerical simulations show that the same asymptotic relaxation holds for the variances $\mathbf{v}(t)$, which converge to $\mathbf{v}^\eq(t)$.

%%%%%%%%%%%%%%%%%%%%%%%%%%%%%
%%%%%%%%%%%%%%%%%%%%%%%%%%%%%
%%%%%%%%%%%%%%%%%%%%%%%%%%%%%
%%%%%%%%%%%%%%%%%%%%%%%%%%%%%

\begin{remark} \label{remark:equilibrium}
Our results on the asymptotic behavior of $\mathbf{m}(t)$ and $\mathbf{v}(t)$ suggest that the solution $\mathbf{f}(\x,t)$ to the Fokker--Planck system \eqref{eq:Fokker-Planck logistic} relaxes over time toward the local equilibrium state $\mathbf{f}^\eq(\x,t)$ defined by \eqref{eq:equilibrium states logistic}. Moreover, the two distributions are expected to converge toward a common global equilibrium state $\mathbf{f}^\infty(\x)$ which is the unique uniform-in-time Gamma distribution canceling out the right-hand sides of system \eqref{eq:Fokker-Planck logistic}, and is given by 
\begin{equation*}
    \begin{split}
        f_1^\infty(x) &= C_1^\infty x^{\frac{2\alpha(\chi+1)}{\sigma_1}\frac{m_1^\infty}{m_1^\infty + m_2^\infty} - 1} \exp\left\{-\frac{2}{\sigma_1 (m_1^\infty + m_2^\infty)} \left(\beta m_2^\infty + \frac{\alpha}{K} m_1^\infty + \alpha \chi \right) x \right\}, \\[4mm]
        f_2^\infty(y) &= C_2^\infty y^{\frac{2 \nu(\theta+1)}{\sigma_2}\frac{m_2^\infty}{m_1^\infty} - 1} \exp\left\{-\frac{2}{\sigma_2 m_1^\infty} (\gamma(\mu - m_1^\infty + \nu\theta) y \right\},
    \end{split}
\end{equation*}
where $C_1^\infty$, $C_2^\infty > 0$ are normalization coefficients depending on the equilibrium means $m_1^\infty$ and $m_2^\infty$. In particular, the means and variances of this equilibrium state give exactly the equilibrium points $\mathbf{m}^\infty$ and $\mathbf{v}^\infty$, namely 
\begin{align*}
    & \int_{\R_+} x f_1^\infty(x) \dd x = m_1^\infty,
    & \int_{\R_+} y f_2^\infty(y) \dd y = m_2^\infty, \qquad \\[4mm]
    \int_{\R_+} & (x - m_1^\infty)^2 f_1^\infty(x) \dd x = v_1^\infty,
    & \int_{\R_+} (y - m_2^\infty)^2 f_2^\infty(y) \dd y = v_2^\infty.
\end{align*}
It remains an open question to determine the equilibration rate of $\mathbf{f}(\x,t)$ toward this global equilibrium state $\mathbf{f}^\infty(\x)$.
\end{remark}

\section*{Conclusion}
In this work we studied the evolution properties of a system of kinetic equations of Fokker--Planck-type, generalizing an analogous model introduced recently in \cite{TosZan} and justifying its derivation from a Boltzmann-type description. These Fokker--Planck equations describe the time evolution of vector-densities accounting for the distribution of the sizes of two populations interacting along the rules of predator--prey-like dynamics. The main feature of the kinetic system is the presence of time-dependent coefficients of drift and diffusion, which are proportional to the mean values of the solution densities, leading to the formation of quasi-equilibria. The importance of the characteristics of these local equilibria in connection with the true solution of the system have been investigated by showing their possible rule. An interesting by-product of the present analysis is the possibility of extending the classical Lotka--Volterra system, which in this case describes the time evolution of the mean solution to a Fokker–Planck system modeling the sizes of predator and prey populations. Specifically, we derive a closed evolution equation for the variances of the solutions, providing additional information on the study of such biological systems. Determining the equilibration rate for the obtained kinetic dynamics is an open and significant question that deserves further attention. In future research, with the aim of quantifying the convergence rate solutions to the kinetic system toward equilibrium, we plan to design suitable entropy methods following the strategy proposed in \cite{ATZ}, where relaxation to equilibrium is shown to occur at least with a polynomial rate. Extensions of the present methods to account for multiple interactions between species, as for example in \cite{BobCerGam,TosTosZan}, are currently under investigation and will be addressed in future works.

%%%%%%%%%%%%%%%%%%%%%%%%%%%%%%%%%%%%%%%%%%%%%%%%%%%%%%%%%%%%%
%%%%%%%%%%%%%%%%%%%%%%%%%%%%%%%%%%%%%%%%%%%%%%%%%%%%%%%%%%%%%
%%%%%%%%%%%%%%%%%     AKNOWLEDGEMENTS    %%%%%%%%%%%%%%%%%%%%
%%%%%%%%%%%%%%%%%%%%%%%%%%%%%%%%%%%%%%%%%%%%%%%%%%%%%%%%%%%%%
%%%%%%%%%%%%%%%%%%%%%%%%%%%%%%%%%%%%%%%%%%%%%%%%%%%%%%%%%%%%%

\bigskip
\noindent \textbf{Acknowledgments.} This work has been written within the activities of GNFM group of INdAM (National Institute of High Mathematics). A.B. acknowledges the support from the European Union’s Horizon Europe research and innovation programme, under the Marie Skłodowska-Curie grant agreement No. 101110920, project MesoCroMo (A Mesoscopic approach to Cross-diffusion Modelling in population dynamics). M.Z. acknowledges partial support by PRIN2022PNRR project No.P2022Z7ZAJ, European Union - NextGenerationEU and by ICSC - Centro Nazionale di Ricerca in High Performance Computing, Big Data and Quantum Computing, funded by European Union - NextGenerationEU.

%%%%%%%%%%%%%%%%%%%%%%%%%%%%%
%%%%%%%%%%%%%%%%%%%%%%%%%%%%%
\begin{figure}[h!]
\begin{flushleft}
\includegraphics[scale=0.3]{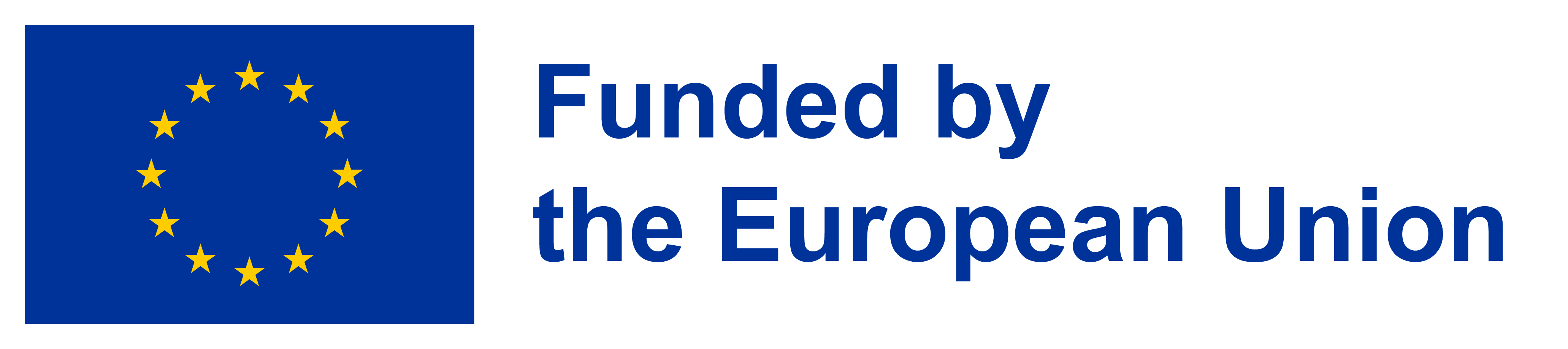}
\end{flushleft}
\end{figure}
%%%%%%%%%%%%%%%%%%%%%%%%%%%%%
%%%%%%%%%%%%%%%%%%%%%%%%%%%%%

\bigskip
\noindent \textbf{Disclaimer.} Funded by the European Union. Views and opinions expressed are however those of the author(s) only and do not necessarily reflect those of the European Union or of the European Research Executive Agency (REA). Neither the European Union nor the granting authority can be held responsible for them.

%%%%%%%%%%%%%%%%%%%%%%%%%%%%%%%%%%%%%%%%%%%%%%%%%%%%%%%%%%%%
%%%%%%%%%%%%%%%%%%%%%%%%%%%%%%%%%%%%%%%%%%%%%%%%%%%%%%%%%%%%
%%%%%%%%%%%%%%%%%%%%    BIBLIOGRAPHY   %%%%%%%%%%%%%%%%%%%%%
%%%%%%%%%%%%%%%%%%%%%%%%%%%%%%%%%%%%%%%%%%%%%%%%%%%%%%%%%%%%
%%%%%%%%%%%%%%%%%%%%%%%%%%%%%%%%%%%%%%%%%%%%%%%%%%%%%%%%%%%%

%\bigskip
%\nocite{*}
\bibliographystyle{plain}
\bibliography{Bibliography_KLV}
\bigskip
\bigskip

\setlength\parindent{0pt}

\end{document}